# POTENTIAL THEORY FOR HYPERBOLIC SPDES[1]

By Robert C. Dalang and Eulalia Nualart

*École Polytechnique Fédérale de Lausanne*

We give general sufficient conditions which imply upper and lower bounds for the probability that a multiparameter process hits a given set $E$ in terms of a capacity of $E$ related to the process. This extends a result of Khoshnevisan and Shi [*Ann. Probab.* **27** (1999) 1135–1159], where estimates for the hitting probabilities of the $(N, d)$ Brownian sheet in terms of the $(d − 2N)$ Newtonian capacity are obtained, and readily applies to a wide class of Gaussian processes. Using Malliavin calculus and, in particular, a result of Kohatsu-Higa [*Probab. Theory Related Fields* **126** (2003) 421–457], we apply these general results to the solution of a system of $d$ nonlinear hyperbolic stochastic partial differential equations with two variables. We show that under standard hypotheses on the coefficients, the hitting probabilities of this solution are bounded above and below by constants times the $(d − 4)$ Newtonian capacity. As a consequence, we characterize polar sets for this process and prove that the Hausdorff dimension of its range is $\min(d, 4)$ a.s.

**1. Introduction.** In this article, we are interested in the following basic problem of potential theory for $\mathbb{R}^d$-valued multiparameter stochastic processes: *given* $E \subset \mathbb{R}^d$, *does this process visit* (*or hit*) $E$ *with positive probability*? Sets that, with probability 1, are not visited are said to be *polar* for the process and otherwise are *nonpolar*. One objective is to relate these hitting probabilities to an analytic expression which is determined by the "geometry" of the set, namely the *capacity* of the set. Another objective is to characterize polar sets for the process. In this article, our main goal is to address these questions for non-Gaussian processes that are solutions to a class of nonlinear hyperbolic stochastic partial differential equations (SPDEs) in the plane (a wide class of Gaussian processes is also considered).

Received September 2002; revised August 2003.

[1]Supported in part by the Swiss National Foundation for Scientific Research.

There is a large literature about potential theory for multiparameter processes. For multiparameter processes whose components are independent single-parameter Markov processes, Fitzsimmons and Salisbury [7] obtained upper and lower bounds on hitting probabilities in terms of a notion of energy of a set. This type of multiparameter process arises in the study of multiple points of single-parameter processes. Song [26] characterized polar sets for the $n$-parameter Ornstein–Uhlenbeck process on a separable Fréchet Gaussian space as null $c_{n,2}$ capacity sets, where the capacity $c_{n,2}$ is defined in a variational form. Hirsch and Song [8] obtained bounds for the hitting probabilities of a class of multiparameter symmetric Markov processes in terms of capacity, also in a variational form.

However, in these references, the class of multiparameter Markov processes does not readily cover certain basic multiparameter processes such as the Brownian sheet or the fractional Brownian sheet. In [10], Khoshnevisan developed a potential theory for a class of multiparameter Markov processes which includes these processes.

This article is essentially motivated by the work of Khoshnevisan and Shi [11], who obtained bounds for hitting probabilities of the Brownian sheet. In particular, if $W = (W_t, t \in \mathbb{R}_+^N)$ denotes an $\mathbb{R}^d$-valued Brownian sheet, they showed that for any compact subset $A$ of $\mathbb{R}^d$ and any $0 < a < b < \infty$, there exists a finite positive constant $K$ such that

$$K^{-1} \operatorname{Cap}_{d-2N}(A) \leq \mathbb{P}\{\exists t \in [a,b]^N : W_t \in A\} \leq K \operatorname{Cap}_{d-2N}(A),$$

where $\operatorname{Cap}_{d-2N}$ denotes the capacity with respect to the Newtonian $(d-2N)$ kernel. The proof of the lower bound is essentially based on estimates of the first and second moments of functionals of occupation measures. The upper bound uses Cairoli's maximal inequality for $N$-parameter martingales as a key step.

In this article, we begin by extending their result to a wide class of $\mathbb{R}^d$-valued continuous multiparameter processes $X = (X_t, t \in \mathbb{R}_+^N)$ that are not necessarily Gaussian but that have absolutely continuous univariate and bivariate distributions away from the axes. In Section 2, we give sufficient conditions on the density of the process that imply upper and lower bounds for the hitting probabilities of $X$ in terms of a given capacity related to the process (Theorem 2.4). For the lower bound on the hitting probability, we require some positivity of a functional of the density of the process (see Hypothesis H1) and an upper bound on a functional of the bivariate density of the process (Hypothesis H2). For the upper bound on the hitting probability, we require that the process be adapted to a commuting filtration (so that Cairoli's maximal inequality can be used) and we need a lower bound for the conditional density of the increment of the process given the past (Hypothesis H3). As a consequence of Theorem 2.4, we obtain an analytic criterion for polarity which is given in Corollary 2.5.

As a first application of the general result of Theorem 2.4, we consider multiparameter Gaussian processes in Section 3. We give sufficient conditions on the covariance function of a Gaussian process which imply bounds for hitting probabilities in terms of the Newtonian capacity (Theorem 3.1). This theorem contains many results that exist in the literature and readily applies to multiparameter Gaussian processes such as the Brownian sheet, $\alpha$-regular Gaussian fields, the Ornstein–Uhlenbeck sheet and the fractional Brownian sheet. For the second and the fourth process, we obtain only a lower bound. The upper bound cannot be obtained from Theorem 3.1 since such processes are not necessarily adapted to a commuting filtration.

In Section 4, we apply results of Malliavin calculus and, in particular, the very recent result of Kohatsu-Higa [13] to the system of nonlinear hyperbolic SPDEs,

$$(1.1) \qquad \begin{aligned} \frac{\partial^2 X_t^i}{\partial t_1 \partial t_2} &= \sum_{j=1}^{d} \sigma_j^i(X_t) \frac{\partial^2 W_t^j}{\partial t_1 \partial t_2} + b^i(X_t), \qquad t = (t_1, t_2) \in \mathbb{R}_+^2, \\ X_t^i &= x_0 \qquad \text{if } t_1 t_2 = 0, 1 \leq i \leq d, \end{aligned}$$

where $W = (W^j, j = 1, \ldots, d)$ is a two-parameter $d$-dimensional Wiener process, the second-order mixed derivative of $W^j$ is the white noise on the plane and $\sigma_j^i$, $b^i$ are smooth functions on $\mathbb{R}^d$. It is known [21] that (1.1) has a unique continuous solution $X = (X_t, t \in \mathbb{R}_+^2)$. In this article, we consider this system of equations in the integral form (4.7) as it was studied in [21]. In the case $b \equiv 0$, under some regularity and strong ellipticity conditions on the matrix $\sigma$ (Conditions P1 and P2), we give in Proposition 4.13 an upper bound of Gaussian type for the conditional density of an increment of $X$ given the past. This uses existing techniques of Malliavin calculus that are adapted to the present context (cf. [17] and [18], Chapter 2). We then use the result of Kohatsu-Higa [13] to establish a Gaussian-type lower bound for the density of the random variable $X_t$ for any $t$ away from the axes (Theorem 4.14) and use a Gaussian-type lower bound for the conditional density of the increment of the process given the past (Theorem 4.16).

In the last section, we apply the results obtained in Sections 2–4 to the solution of system (1.1). In the case $b \equiv 0$, we prove in Theorem 5.1 that under Conditions P1 and P2 introduced in Section 4.4, the hitting probabilities of the solution can be bounded above and below in terms of the $(d-4)$ Newtonian capacity. The verification of Hypothesis H2 uses the upper bound of Gaussian type obtained in Section 4.4. The main effort in proving Theorem 5.1 lies in verification of Hypothesis H3, which uses the lower bounds of Gaussian type for the density of the solution obtained in Section 4.5. These Gaussian-type lower bounds also imply the positivity of the density of the solution and so the verification of Hypothesis H1. We treat the case $b \not\equiv 0$ via a change of probability measure (see Corollary 5.3). As a consequence of Corollary 5.3, we prove in Corollaries 5.4 and 5.5 that polar sets for the solution to (1.1) are those of $(d-4)$ Newtonian capacity zero, that the Hausdorff dimension of the range of the solution is $\min(d,4)$ almost surely and its stochastic codimension is $(d-4)^+$. Finally, we identify $d=3$ as the critical dimension for the solution to hit points in $\mathbb{R}^d$ (see Corollary 5.6).

Notice that we obtain the same zero capacity condition for polarity and Hausdorff dimension obtained by Dynkin [6], LeGall [16] and Perkins [24] for super-Brownian motion. However, there is no direct connection between this work and theirs.

## 2. General theory.

For any $s,t \in \mathbb{R}_+^N$, we write $s \le t$ when $s_i \le t_i$ for all $i=1,\ldots,N$, where $s_i$ denotes the $i$th coordinate of $s$, and we write $s < t$ when $s \le t$ and $s \neq t$. By $s \wedge t$, we mean the point whose $i$th coordinate is $s_i \wedge t_i$ for all $i=1,\ldots,N$. If $s,t \in \mathbb{R}_+^N$ with $s \le t$, we write $[s,t] = \prod_{i=1}^N [s_i,t_i]$ for an $N$-dimensional rectangle and $(s,t] = \prod_{i=1}^N (s_i,t_i]$. Finally, we denote by $\|\cdot\|$ the Euclidean norm.

Let $(\Omega,\mathcal{G},\mathbb{P})$ be a complete probability space and let $\mathcal{F} = (\mathcal{F}_t, t \in \mathbb{R}_+^N)$ be a complete, right continuous, commuting filtration, that is, an increasing family of sub-$\sigma$-fields of $\mathcal{G}$ such that:

(i) $\mathcal{F}_0$ contains all the null sets of $\mathcal{G}$.

(ii) For every $t \in \mathbb{R}_+^N$, $\mathcal{F}_t = \bigcap_{s>t} \mathcal{F}_s$.

(iii) For every $s,t \in \mathbb{R}_+^N$ and for all bounded, $\mathcal{F}_t$-measurable random variables $Y$,

$$\mathbb{E}[Y|\mathcal{F}_s] = \mathbb{E}[Y|\mathcal{F}_{s \wedge t}] \qquad \text{a.s.}$$

Note that when $N=2$, (iii) is hypothesis (F4) of Cairoli and Walsh [1]. For $N>2$, hypothesis (iii) appears in [10], Chapter 7, Section 2.1.

Let $X = (X_t, t \in \mathbb{R}_+^N)$ be a continuous $\mathbb{R}^d$-valued stochastic process defined on $(\Omega,\mathcal{G},\mathbb{P})$ and *not necessarily* adapted to $\mathcal{F}$. We suppose that for all $s,t \in (0,+\infty)^N$ with $t \neq s$, the distribution of the random variable $(X_t, X_s)$ has a density that is denoted $p_{X_t,X_s}(x,y)$. We write $p_{X_t}(x)$ for the density of the random variable $X_t$ for all $t \in (0,+\infty)^N$.

Given a Borel subset $E$ of $\mathbb{R}^d$, we denote by $\mathcal{P}(E)$ the collection of all probability measures on $\mathbb{R}^d$ with support in $E$.

Following [10], Appendix D, we say that $k(\cdot)$ is a *kernel* in $\mathbb{R}^d$ (or a *gauge function*) if $k(\cdot)$ is an even, nonnegative and locally integrable function on $\mathbb{R}^d$ which is continuous on $\mathbb{R}^d \setminus \{0\}$ and positive in a neighborhood of the origin. Basic examples of kernels are the Newtonian $\beta$ kernels $k_\beta(\cdot)$ (see Section 3.1).

Given a kernel $k(\cdot)$, for any $\mu \in \mathcal{P}(E)$, we write

$$\mathcal{E}_k(\mu) = \int_{\mathbb{R}^d} \int_{\mathbb{R}^d} k(x-y)\mu(dx)\mu(dy)$$

and term this quantity the $k$ *energy* of $\mu$. The $k$ *capacity* of $E$ is defined by

$$\text{Cap}_k(E) = \frac{1}{\inf_{\mu \in \mathcal{P}(E)} \mathcal{E}_k(\mu)}.$$

The following properties of $\text{Cap}_k(\cdot)$ are given in [10], Appendix D, in particular, Lemma 2.1.2 there.

LEMMA 2.1. $\text{Cap}_k(\cdot)$ *has the following properties*:

(a) Monotonicity. *For any two Borel subsets $E_1 \subset E_2$ of $\mathbb{R}^d$, $\text{Cap}_k(E_1) \leq \text{Cap}_k(E_2)$.*

(b) Outer regularity on compact sets. *For any sequence $E, E_1, E_2, \ldots$ of compact subsets of $\mathbb{R}^d$ such that $E_n \downarrow E$, $\lim_{n\to\infty} \text{Cap}_k(E_n) = \text{Cap}_k(E)$.*

Let $\mathcal{P}_0(E)$ denote the collection of all probability measures on $\mathbb{R}^d$ with support in $E$ that are absolutely continuous with respect to Lebesgue measure. The *absolutely continuous capacity* $\text{Cap}_k^0(E)$ of $E$ with respect to $k(\cdot)$ is defined by

$$\text{Cap}_k^0(E) = \frac{1}{\inf_{\mu \in \mathcal{P}_0(E)} \mathcal{E}_k(\mu)}.$$

Since $\text{Cap}_k^0$ is not outer regular on compact sets (cf. [10], Appendix D, Section 2.2), we define, for all bounded Borel sets $E \subset \mathbb{R}^d$,

$$\text{Cap}_k^{\text{ac}}(E) = \inf\{\text{Cap}_k^0(\overline{F}) : F \supset E, F \text{ bounded and open}\}.$$

Then $\text{Cap}_k(A) \geq \text{Cap}_k^{\text{ac}}(A)$. We now state an additional condition on $k$ (which is related to the classical notion of *balayage*; see [15], Chapter IV) that ensures that capacity and absolutely continuous capacity with respect to $k$ agree on compact sets.

Following [10], Appendix D, we say that a kernel $k(\cdot)$ on $\mathbb{R}^d$ is *proper* if, for all compact sets $A \subset \mathbb{R}^d$ and $\mu \in \mathcal{P}(A)$, there exist bounded open sets $A_1, A_2, \ldots$ such that:

1. $\overline{A}_n \downarrow A$.
2. For all large $n \geq 1$, there exist absolutely continuous measures $\mu_n$ with support contained in $\overline{A}_n$ such that, for all $\varepsilon > 0$, there exists $N_0$ such that for all $n \geq N_0$:

   (a) $\mu_n(\overline{A}_n) \geq 1 - \varepsilon$;
   (b) $\int_{\mathbb{R}^d} k(x - y) \mu_n(dy) \leq \int_{\mathbb{R}^d} k(x - y) \mu(dy)$ for all $x \in \mathbb{R}^d$.

PROPOSITION 2.2 ([10], Appendix D, Theorem 2.3.1). *Let $k(\cdot)$ be a proper kernel in $\mathbb{R}^d$. Then, for all compact sets $A \subset \mathbb{R}^d$, $\text{Cap}_k(A) = \text{Cap}_k^{\text{ac}}(A)$.*

We now introduce the following hypotheses, which ensure a lower bound on hitting probabilities for $X$ [see Theorem 2.4(a)].

HYPOTHESIS H1.  *For all $0 < a < b < \infty$ and $M > 0$, there exists a finite positive constant $C_1(a, b, M)$ such that for almost all $x \in \mathbb{R}^d$ with $\|x\| \leq M$,*

$$\int_{[a,b]^N} p_{X_t}(x) \, dt \geq C_1.$$

HYPOTHESIS H2.  *There exists a proper kernel $k(\cdot)$ in $\mathbb{R}^d$ such that for all $0 < a < b < \infty$ and $M > 0$, there exists a finite positive constant $C_2(a, b, M)$ such that, for almost all $x, y \in \mathbb{R}^d$ with $\|x\| \leq M$ and $\|y\| \leq M$,*

$$\int_{[a,b]^N} \int_{[a,b]^N} p_{X_t, X_s}(x, y) \, dt \, ds \leq C_2 k(x - y).$$

In the case where $X$ is adapted to $\mathcal{F}$, for $s \in (0, +\infty)^N$, let $P_s(\omega, \cdot)$ be a regular version of the conditional distribution of the process $(X_t - X_s, t \in \mathbb{R}_+^N \setminus [0, s])$ given $\mathcal{F}_s$. If for almost all $\omega$ and all $t \in \mathbb{R}_+^N \setminus [0, s]$, the law of $X_t - X_s$ under $P_s(\omega, \cdot)$ is absolutely continuous with respect to Lebesgue measure on $\mathbb{R}^d$, we let $p_{s,t}(\omega, x)$ denote the density of $X_t - X_s$ under $P_s(\omega, \cdot)$. In this case, there is a null set $N_s \in \mathcal{F}_s$ such that for $\omega \in \Omega \setminus N_s$, $E$ a Borel subset of $\mathbb{R}^d$ and $s < t$,

$$P_s(\omega, \{X_t - X_s \in E\}) = \int_E p_{s,t}(\omega, x) \, dx.$$

In particular, $p_{s,t}(\omega, x)$ is a version of the conditional density of $X_t - X_s$ given $\mathcal{F}_s$. The function $(\omega, t, x) \mapsto p_{s,t}(\omega, x)$ can be chosen to be measurable.

PROPOSITION 2.3.  *Let $f : \mathbb{R}^d \times \mathbb{R}^d \to \mathbb{R}$ be a nonnegative Borel function, let $Y$ be an $\mathcal{F}_s$-measurable random variable and suppose that $\mathbb{E}[f(X_t - X_s, Y)] < \infty$. Then*

$$\mathbb{E}[f(X_t - X_s, Y) | \mathcal{F}_s] = \int_{\mathbb{R}^d} f(x, Y) p_{s,t}(\omega, x) \, dx \qquad a.s.$$

PROOF.  If $f(x, y) = f_1(x) f_2(y)$, then using [5], Theorem 10.2.5, we have

$$\begin{aligned}
\mathbb{E}[f(X_t - X_s, Y) | \mathcal{F}_s] &= \mathbb{E}[f_1(X_t - X_s) f_2(Y) | \mathcal{F}_s] \\
&= f_2(Y) \mathbb{E}[f_1(X_t - X_s) | \mathcal{F}_s] \\
&= f_2(Y) \int_{\mathbb{R}^d} f_1(x) p_{s,t}(\omega, x) \, dx \qquad \text{a.s.} \\
&= \int_{\mathbb{R}^d} f(x, Y) p_{s,t}(\omega, x) \, dx \qquad \text{a.s.}
\end{aligned}$$

One easily concludes the proof using a monotone class argument ([3], Chapter I, Theorem 21). □

We now introduce a third hypothesis, which leads to an upper bound on hitting probabilities for $X$ [see Theorem 2.4(b)].

HYPOTHESIS H3. *For all $0 < a < b < \infty$ and $M > 0$, there exists a finite positive constant $C_3(a, b, M)$ such that for all $s \in [a, b]^N$, a.s., for almost all $x \in \mathbb{R}^d$,*

$$\int_{[b, 2b-a]^N} p_{s,t}(\omega, x) \, dt \geq C_3 k(x) \mathbb{1}_{\{\|x + X_s\| \leq M, \|X_s\| \leq M\}}(\omega),$$

*where $k(x)$ is the same kernel as in Hypothesis H2.*

We are now ready to state the main result of this section.

THEOREM 2.4. (a) *Assuming Hypotheses H1 and H2, for all $0 < a < b < \infty$ and $M > 0$, there exists a finite positive constant $K_1(a, b, M)$ such that for all compact sets $A \subset \{x \in \mathbb{R}^d : \|x\| < M\}$,*

$$K_1 \operatorname{Cap}_k(A) \leq \mathbb{P}\{\exists t \in [a, b]^N : X_t \in A\}.$$

(b) *If $(X_t, t \in \mathbb{R}_+^N)$ is adapted to a commuting filtration $(\mathcal{F}_t, t \in \mathbb{R}_+^N)$, and Hypotheses H2 and H3 hold, then for all $0 < a < b < \infty$ and $M > 0$, there exists a finite positive constant $K_2(a, b, M)$ such that for all compact sets $A \subset \{x \in \mathbb{R}^d : \|x\| < M\}$,*

$$\mathbb{P}\{\exists t \in [a, b]^N : X_t \in A\} \leq K_2 \operatorname{Cap}_k(A).$$

Before proving Theorem 2.4, we mention an important consequence. Recall that a Borel set $E \subset \mathbb{R}^d$ is said to be *polar* for the process $X$ if

$$\mathbb{P}\{\exists t \in (0, +\infty)^N : X_t \in E\} = 0.$$

COROLLARY 2.5. *For a process $X$ adapted to a commuting filtration, under Hypotheses H1–H3, a compact subset $E$ of $\mathbb{R}^d$ is polar for $X$ if and only if $\operatorname{Cap}_k(E) = 0$.*

PROOF. If $E$ is polar for $X$, then clearly $\operatorname{Cap}_k(E) = 0$ by Theorem 2.4(a). Conversely, suppose $\operatorname{Cap}_k(E) = 0$. Write $(0, +\infty)^N = \bigcup_{m \in \mathbb{N}} [\frac{1}{m}, m]^N$. By Theorem 2.4(b), for all $m \geq 1$, there is $K_2 < \infty$ (depending on $m$) such that

$$\mathbb{P}\left\{\exists t \in \left[\frac{1}{m}, m\right]^N : X_t \in E\right\} \leq K_2 \operatorname{Cap}_k(E) = 0.$$

Since this holds for all $m$, $E$ is polar for $X$.   $\square$

PROOF OF THEOREM 2.4.   (a) *The lower bound.* Suppose $0 < a < b < \infty$ and $0 < M < \infty$ are fixed. Let $A$ be a compact subset of $\{x \in \mathbb{R}^d : \|x\| < M\}$. For $\varepsilon \in (0,1)$, define $\overline{A}_\varepsilon = \{x \in \mathbb{R}^d : \mathrm{dist}(x, A) \le \varepsilon\}$, the closed $\varepsilon$ enlargement of $A$, where $\mathrm{dist}(x, A) = \|x - \mathrm{proj}_A x\|$ and $\mathrm{proj}_A x$ denotes the orthogonal projection of $x$ on $A$. Fix $\varepsilon \in (0,1)$ and let $f$ be a probability density on $\mathbb{R}^d$ whose support is contained in $\overline{A}_\varepsilon$. We consider the functional $J_{a,b}(f)$ defined by

$$J_{a,b}(f) = \int_{[a,b]^N} f(X_t) \, dt.$$

By the Cauchy–Schwarz inequality,

(2.1)     $\mathbb{P}\{\exists \, t \in [a,b]^N : X_t \in \overline{A}_\varepsilon\} \ge \mathbb{P}\{J_{a,b}(f) > 0\} \ge \dfrac{\{\mathbb{E}[J_{a,b}(f)]\}^2}{\mathbb{E}[\{J_{a,b}(f)\}^2]}.$

Using Fubini's theorem, we easily deduce from Hypothesis H1 that

(2.2)                           $\mathbb{E}[J_{a,b}(f)] \ge C_1$

and from Hypothesis H2 that

(2.3)                       $\mathbb{E}[\{J_{a,b}(f)\}^2] \le C_2 \mathcal{E}_k(f),$

where $\mathcal{E}_k(f)$ denotes the $k$ energy of the measure $f(x) \, dx$. Applying (2.2) and (2.3) to (2.1), we obtain

$$\mathbb{P}\{\exists \, t \in [a,b]^N : X_t \in \overline{A}_\varepsilon\} \ge \frac{C_1^2}{C_2 \mathcal{E}_k(f)}.$$

Take the supremum over all $f(x) \, dx \in \mathcal{P}(\overline{A}_\varepsilon)$ and see that for all $\varepsilon > 0$,

$$\mathbb{P}\{\exists \, t \in [a,b]^N : X_t \in \overline{A}_\varepsilon\} \ge \frac{C_1^2}{C_2} \mathrm{Cap}_k^{\mathrm{ac}}(\overline{A}_\varepsilon).$$

By Proposition 2.2, we can replace $\mathrm{Cap}_k^{\mathrm{ac}}(\overline{A}_\varepsilon)$ by $\mathrm{Cap}_k(\overline{A}_\varepsilon)$ because $k$ is proper. As $\varepsilon \to 0^+$, $\overline{A}_\varepsilon \downarrow A$, which is compact. By Lemma 2.1(b), $\mathrm{Cap}_k(\overline{A}_\varepsilon)$ converges to $\mathrm{Cap}_k(A)$ as $\varepsilon \to 0^+$. Finally, since $A$ is compact and the process $t \mapsto X_t$ is continuous,

$$\bigcap_{\varepsilon > 0} \{\exists \, t \in [a,b]^N : X_t \in \overline{A}_\varepsilon\} = \{\exists \, t \in [a,b]^N : X_t \in A\}.$$

We conclude that

$$\mathbb{P}\{\exists \, t \in [a,b]^N : X_t \in A\} \ge \frac{C_1^2}{C_2} \mathrm{Cap}_k(A).$$

This concludes the proof of (a) of Theorem 2.4.

(b) *The upper bound.* Suppose $0 < a < b < \infty$ and $0 < M < \infty$ are fixed. Let $A$ be a compact subset of $\{x \in \mathbb{R}^d : \|x\| < M\}$. Let $f : \mathbb{R}^d \to \mathbb{R}_+$ be a measurable function such that $\mathbb{E}[\{J_{b,2b-a}(f)\}^2] < \infty$. We define the following square integrable multiparameter martingale:

$$M_t(f) = \mathbb{E}[J_{b,2b-a}(f)|\mathcal{F}_t], \qquad t \in \mathbb{R}_+^N.$$

Since $\mathcal{F}_t$ is a commuting filtration, we can apply Cairoli's maximal inequality (see [10], Chapter 7, Theorem 2.3.2), to get

$$\mathbb{E}\left[\sup_{t \in [a,b]^N \cap \mathbb{D}^N} \{M_t(f)\}^2\right] \le 4^N \sup_{t \in [a,b]^N \cap \mathbb{D}^N} \mathbb{E}[\{M_t(f)\}^2] \le 4^N \mathbb{E}[\{J_{b,2b-a}(f)\}^2],$$

where $\mathbb{D}$ denotes the set of dyadic rationals. Suppose that $f$ is a density function on $\mathbb{R}^d$ supported on $\{x \in \mathbb{R}^d : \|x\| < M\}$. By Hypothesis H2 and (2.3), we get

$$(2.4) \qquad \mathbb{E}\left[\sup_{t \in [a,b]^N \cap \mathbb{D}^N} \{M_t(f)\}^2\right] \le 4^N C_2 \mathcal{E}_k(f).$$

For $\varepsilon \in (0,1)$, define $A_\varepsilon = \{x \in \mathbb{R}^d : \mathrm{dist}(x, A) < \varepsilon\}$, the open $\varepsilon$ enlargement of $A$. Suppose $\varepsilon$ is small enough so that $A_\varepsilon \subset \{x \in \mathbb{R}^d : \|x\| < M\}$. We can assume that

$$(2.5) \qquad \mathbb{P}\{\exists\, t \in [a,b]^N : X_t \in A_\varepsilon\} > 0 \qquad \text{for all } \varepsilon > 0.$$

Indeed, if there exists an $\varepsilon > 0$ such that this probability is equal to zero, then the upper bound is trivial since

$$\mathbb{P}\{\exists\, t \in [a,b]^N : X_t \in A\} \le \mathbb{P}\{\exists\, t \in [a,b]^N : X_t \in A_\varepsilon\} \qquad \text{for all } \varepsilon > 0.$$

Assuming (2.5), we claim that there exists a random variable $T^\varepsilon$, taking values in $([a,b]^N \cap \mathbb{D}^N) \cup \{+\infty\}$, such that

$$(2.6) \qquad \{T^\varepsilon < \infty\} \quad \Longleftrightarrow \quad \{\exists\, t \in [a,b]^N \cap \mathbb{D}^N : X_t \in A_\varepsilon\}$$

and $X_{T^\varepsilon} \in A_\varepsilon$ on $\{T^\varepsilon < \infty\}$. Indeed, order the set $[a,b]^N \cap \mathbb{D}^N = \{q_1, q_2, \ldots\}$ and define $T^\varepsilon = q_{\inf\{k : X_{q_k} \in A_\varepsilon\}}$, where $\inf \varnothing$ is defined to be $+\infty$ and in this case $T^\varepsilon = +\infty$. Note that assumption (2.5), the fact that $A_\varepsilon$ is open and the continuity of $t \mapsto X_t$ imply that $\mathbb{P}\{T^\varepsilon < \infty\} > 0$. In particular, it is possible to condition on the event $\{T^\varepsilon < \infty\}$.

For any Borel set $E \subset \mathbb{R}^d$, define

$$\mu_\varepsilon(E) = \mathbb{P}\{X_{T^\varepsilon} \in E | T^\varepsilon < \infty\}.$$

Clearly $\mu_\varepsilon \in \mathcal{P}(A_\varepsilon)$. Moreover, $\mu_\varepsilon$ is absolutely continuous because every $X_t$ is: Let $f_\varepsilon(x)$ be the Radon–Nikodym derivative of $\mu_\varepsilon$ with respect to $\lambda$. Then $f_\varepsilon$ is supported on $A_\varepsilon$. Applying (2.4) to $f_\varepsilon$, we have that for all $\varepsilon \in (0,1)$,

$$\mathbb{E}\left[\sup_{t \in [a,b]^N \cap \mathbb{D}^N} \{M_t(f_\varepsilon)\}^2\right] \le 4^N C_2 \mathcal{E}_k(f_\varepsilon).$$

We claim that for all $\varepsilon \in (0,1)$ and all $s \in [a,b]^N$,

$$(2.7) \qquad M_s(f_\varepsilon) \geq C_3 \mathbb{1}_{\{\|X_s\| \leq M\}} \int_{\mathbb{R}^d} f_\varepsilon(x + X_s) k(x) \, dx \qquad \text{a.s.}$$

Indeed, fix $s \in [a,b]^N$. Using the fact that $X_s$ is $\mathcal{F}_s$-measurable, Proposition 2.3 and Hypothesis H3, we have that for almost all $\omega \in \Omega$,

$$M_s(f_\varepsilon) = \mathbb{E}\left[\int_{[b,2b-a]} f_\varepsilon(X_t - X_s + X_s) \, dt \Big| \mathcal{F}_s\right]$$

$$= \int_{\mathbb{R}^d} f_\varepsilon(x + X_s)\left(\int_{[b,2b-a]} p_{s,t}(\omega, x) \, dt\right) dx \qquad \text{a.s.}$$

$$\geq C_3 \mathbb{1}_{\{\|X_s\| \leq M\}} \int_{\mathbb{R}^d} f_\varepsilon(x + X_s) k(x) \, dx \qquad \text{a.s.}$$

The last inequality follows since $f_\varepsilon$ is nonnegative and supported in $A_\varepsilon$. This completes the proof of (2.7).

Since (2.7) holds for all $s \in [a,b]^N$ and since $T^\varepsilon \in [a,b]^N$, we can replace $s$ by $T^\varepsilon$ in (2.7). Note that $\{\|X_{T^\varepsilon}\| \leq M\}$ holds on $\{T^\varepsilon < \infty\}$. Therefore,

$$(2.8) \qquad \sup_{s \in [a,b]^N \cap \mathbb{D}^N} M_s(f_\varepsilon) \geq C_3 \mathbb{1}_{\{T^\varepsilon < \infty\}} \int_{\mathbb{R}^d} f_\varepsilon(x + X_{T^\varepsilon}) k(x) \, dx \qquad \text{a.s.}$$

Square both sides of the last inequality, take expectations and apply (2.4) to the left-hand side, to obtain

$$4^N C_2 \mathcal{E}_k(f_\varepsilon) \geq C_3^2 \mathbb{P}\{T^\varepsilon < \infty\} \mathbb{E}\left[\left(\int_{\mathbb{R}^d} f_\varepsilon(x + X_{T^\varepsilon}) k(x) \, dx\right)^2 \Big| T^\varepsilon < \infty\right]$$

$$= C_3^2 \mathbb{P}\{T^\varepsilon < \infty\} \int_{\mathbb{R}^d}\left(\int_{\mathbb{R}^d} f_\varepsilon(x + y) k(x) \, dx\right)^2 f_\varepsilon(y) \, dy.$$

Using Jensen's inequality, we get

$$4^N C_2 \mathcal{E}_k(f_\varepsilon) \geq C_3^2 \mathbb{P}\{T^\varepsilon < \infty\}\left(\int_{\mathbb{R}^d} \int_{\mathbb{R}^d} k(x - y) f_\varepsilon(x) f_\varepsilon(y) \, dx \, dy\right)^2$$

$$= C_3^2 \mathbb{P}\{T^\varepsilon < \infty\}(\mathcal{E}_k(f_\varepsilon))^2.$$

If $\mathcal{E}_k(f_\varepsilon)$ were finite, this would imply

$$(2.9) \qquad \mathbb{P}\{T^\varepsilon < \infty\} \leq \frac{4^N C_2}{C_3^2 \mathcal{E}_k(f_\varepsilon)},$$

but we do not know a priori that $f_\varepsilon$ has finite energy. For that reason, we use a truncation argument.

For all $q > 0$ and all $\varepsilon \in (0,1)$, define

$$f_\varepsilon^q(x) = f_\varepsilon(x) \mathbb{1}_{[0,q]}(f_\varepsilon(x)), \qquad x \in \mathbb{R}^d.$$

Since $f_\varepsilon$ is supported on $A_\varepsilon$, so is $f_\varepsilon^q$. Moreover, the latter is a subprobability density function that is bounded above by $f_\varepsilon$ and $q$. Therefore, since $k$ is locally integrable in $\mathbb{R}^d$, $\mathcal{E}_k(f_\varepsilon^q) < \infty$. Apply to $f_\varepsilon^q$ exactly the same argument that led to (2.8) to see that

$$\sup_{s \in [a,b]^N \cap \mathbb{D}^N} M_s(f_\varepsilon^q) \geq C_3 \mathbb{1}_{\{T^\varepsilon < \infty\}} \int_{\mathbb{R}^d} f_\varepsilon^q(x + X_{T^\varepsilon}) k(x)\, dx \qquad \text{a.s.}$$

Square both sides of the inequality, take expectations and use Jensen's inequality to get

$$\mathbb{E}\bigg[ \sup_{s \in [a,b]^N \cap \mathbb{D}^N} \{M_s(f_\varepsilon^q)\}^2 \bigg]$$

$$\geq C_3^2 \mathbb{P}\{T^\varepsilon < \infty\} \bigg( \int_{\mathbb{R}^d} \int_{\mathbb{R}^d} k(x-y) f_\varepsilon^q(x) f_\varepsilon(y)\, dx\, dy \bigg)^2.$$

By (2.4), the left-hand side is bounded above by $4^N C_2 \mathcal{E}_k(f_\varepsilon^q)$. The right-hand side is clearly bounded below by $C_3^2 \mathbb{P}\{T^\varepsilon < \infty\}(\mathcal{E}_k(f_\varepsilon^q))^2$. Hence, we obtain

$$\mathbb{P}\{T^\varepsilon < \infty\} \leq \frac{4^N C_2}{C_3^2 \mathcal{E}_k(f_\varepsilon^q)}.$$

Finally, since $k(x)$ is nonnegative, $\lim_{q \uparrow +\infty} \mathcal{E}_k(f_\varepsilon^q) = \mathcal{E}_k(f_\varepsilon)$, we can let $q \uparrow +\infty$ in the above inequality and use the monotone convergence theorem to obtain (2.9).

Now, since $t \mapsto X_t$ is continuous and $A_\varepsilon$ is open, using (2.6) and (2.9), we obtain

$$\mathbb{P}\{\exists t \in [a,b]^N : X_t \in A_\varepsilon\} \leq \frac{4^N C_2}{C_3^2 \mathcal{E}_k(\mu_\varepsilon)}.$$

Recall that $\mu_\varepsilon \in \mathcal{P}(A_\varepsilon)$, so for all $\varepsilon > 0$,

$$\mathbb{P}\{\exists t \in [a,b]^N : X_t \in A\} \leq \mathbb{P}\{\exists t \in [a,b]^N : X_t \in A_\varepsilon\} \leq \frac{4^N C_2}{C_3^2} \operatorname{Cap}_k(\overline{A}_\varepsilon).$$

Finally, since $\overline{A}_\varepsilon$ is compact, Lemma 2.1(b) implies that $\operatorname{Cap}_k(\overline{A}_\varepsilon)$ converges as $\varepsilon \to 0^+$ to $\operatorname{Cap}_k(A)$. This concludes the proof of Theorem 2.4(b). $\square$

**3. Multiparameter Gaussian processes.** In this section, we focus on Gaussian processes and reformulate the Hypotheses H1–H3 as conditions on the covariance of the process so as to relate bounds on hitting probabilities directly to properties of the covariance.

3.1. *Relating Newtonian capacity and covariance.* Let $(\Omega, \mathcal{G}, \mathbb{P})$ be a complete probability space. Let $X = (X_t, t \in \mathbb{R}_+^N)$ be a continuous $\mathbb{R}^d$-valued centered Gaussian process with independent coordinate processes $(X_t^i)$. For all $t \in (0, \infty)^N$, we denote by $p_{X_t}(x)$ the density function of the centered Gaussian random variable $X_t$ on $\mathbb{R}^d$. For all $s, t \in (0, \infty)^N$, we write $\sigma(s, t) = \mathbb{E}[X_s^i X_t^i]$, which does not depend on $i$, $\sigma^2(t) = \sigma(t, t)$ and $\rho(s, t) = \sigma(s, t)/(\sigma(t)\sigma(s))$.

Given $\alpha \in (0, 1)$ and $\gamma \geq \alpha$, we introduce the following hypotheses.

HYPOTHESIS A1. *For all $0 < a < b < \infty$, there exist positive finite constants $\delta$, $\varepsilon$ and $C_1, \ldots, C_5$ such that for all $s, t \in [a, b]^N$,*

$$(3.1) \qquad C_1 \leq \sigma^2(t) \leq C_2,$$

$$(3.2) \qquad \left| 1 - \frac{\sigma(s, t)}{\sigma^2(s)} \right| \leq C_3 \|t - s\|^\gamma \qquad \text{if } \|t - s\| \leq \delta,$$

$$(3.3) \qquad C_4 \|t - s\|^{2\alpha} \leq 1 - \rho^2(s, t) \leq C_5 \|t - s\|^{2\alpha} \qquad \text{if } \|t - s\| \leq \delta,$$

$$(3.4) \qquad |\rho(s, t)| < 1 - \varepsilon \qquad \text{if } \|t - s\| > \delta.$$

HYPOTHESIS A2. *For all $0 < a < b < \infty$ and $M > 0$, there exist finite positive constants $C_6$, $C_7$ and $C_8$ such that for all $s, t \in [a, b]^N$ with $s \leq t$,*

$$(3.5) \qquad \mathbb{1}_{\{\|X_s\| \leq M\}} |\mathbb{E}[X_t - X_s | \mathcal{F}_s]| \leq C_6 \|t - s\|^\gamma,$$

$$(3.6) \qquad C_7 \|t - s\|^{2\alpha} \leq \mathbb{E}[(X_t - X_s) - \mathbb{E}[X_t - X_s | \mathcal{F}_s]]^2 \leq C_8 \|t - s\|^{2\alpha}.$$

THEOREM 3.1. (a) *Assume there are $\alpha \in (0, 1)$ and $\gamma \geq \alpha$ for which Hypothesis A1 holds and $N/\alpha \geq 2$. Then for all $0 < a < b < \infty$ and $M > 0$, there exists a finite positive constant $K_1(a, b, M)$, such that for all compact sets $A \subset \{x \in \mathbb{R}^d : \|x\| < M\}$,*

$$K_1 \operatorname{Cap}_{d - (N/\alpha)}(A) \leq \mathbb{P}\{\exists t \in [a, b]^N : X_t \in A\},$$

*where, for $\beta \geq 0$, $\operatorname{Cap}_\beta(\cdot)$ denotes the capacity with respect to the Newtonian $\beta$ kernel $k_\beta(\cdot)$, where*

$$k_\beta(x) = \begin{cases} \|x\|^{-\beta}, & \text{if } 0 < \beta < d, \\ \ln\left(\dfrac{3M}{\|x\|}\right), & \text{if } \beta = 0, \end{cases}$$

*and for $\beta < 0$, $\operatorname{Cap}_\beta(\cdot) = 1$.*

(b) *Suppose $(X_t, t \in \mathbb{R}_+^N)$ is adapted to a commuting filtration $(\mathcal{F}_t, t \in \mathbb{R}_+^N)$. Assume there are $\alpha \in (0, 1)$ and $\gamma \geq \alpha$ for which Hypotheses A1 and A2 hold and $N/\alpha \geq 2$. Then for all $0 < a < b < \infty$ and $M > 0$, there exists a finite positive constant $K_2(a, b, M)$ such that for all compact sets $A \subset \{x \in \mathbb{R}^d : \|x\| < M\}$,*

$$\mathbb{P}\{\exists t \in [a, b]^N : X_t \in A\} \leq K_2 \operatorname{Cap}_{d - (N/\alpha)}(A).$$

REMARK 3.2. (a) For $0 < \beta < d$, the functions $k_\beta$ are not only positive but even positive-definite: this follows from [27], Chapter V, Section 1, Lemma 2(b). For $\beta = 0$, the kernel $k_0$ is not nonnegative on $\mathbb{R}^d$, but $k_0(x - y) > 0$ for $x, y \in A$ when $A \subset \{x \in \mathbb{R}^d : \|x\| < M\}$, and this is sufficient for the results of Section 2 to hold.

(b) Note that for any $0 \leq \beta \leq d-2$, $k_\beta(\cdot)$ is a proper kernel. Indeed, given a compact set $A \subset \{x \in \mathbb{R}^d : \|x\| < M\}$ and $\mu \in \mathcal{P}(A)$, for any $n \geq 1$, let $A_n = \{x \in \mathbb{R}^d : \text{dist}(x, A) < 1/n\}$ be the open enlargement of $A$. Let $(B_t, \ t \geq 0)$ be a standard Brownian motion in $\mathbb{R}^d$ and let $p_t$ be the density of $B_t$. For $n \geq 1$, let $\mu_{n,t}$ be the measure whose density function is the restriction to the set $\bar{A}_n$ of $p_t * \mu$, where $*$ denotes the convolution product. Set $f(x) = (\mu * k_\beta)(x)$. Since $\Delta k_\beta(x) \leq 0$ for $0 \leq \beta \leq d-2$, $E_x(f(B_t)) \leq f(x)$ for all $t \geq 0$ and $x \in \mathbb{R}^d$, or equivalently,

$$\int_{\mathbb{R}^d} k_\beta(x - y)(p_t * \mu)(y) \, dy \leq \int_{\mathbb{R}^d} k_\beta(x - y) \, \mu(dy).$$

For small $t$, the $\mu_{n,t}$ are nearly probability measures, and so (a) and (b) of the definition of a proper kernel hold.

(c) The condition $N/\alpha \geq 2$ ensures that $k_{d-(N/\alpha)}$ is a proper kernel. This is only a restriction when $N = 1$.

(d) If $\alpha = 1/2$ and $d < 2N$, the choice $\text{Cap}_{d-2N}(A) = 1$ is natural, since in this case, the Brownian sheet hits points in $\mathbb{R}^d$ (cf. [23]).

Before proving Theorem 3.1, we mention an important consequence. Given $s \geq 0$ and a Borel subset $E$ of $\mathbb{R}^d$, let

$$\mathcal{H}_s(E) = \lim_{\varepsilon \to 0^+} \inf \left\{ \sum_{i=1}^\infty (2r_i)^s : E \subset \bigcup_{i=1}^\infty \mathcal{B}(x_i, r_i), \sup_i r_i \leq \varepsilon \right\},$$

where $\mathcal{B}(x, r)$ denotes the closed ball of radius $r > 0$ centered at $x \in \mathbb{R}^d$. The $\mathcal{H}_s$ is called the *d-dimensional Hausdorff measure*. Moreover, we associate to the set $E$ a number $\dim_{\mathcal{H}}(E)$ as follows:

$$\dim_{\mathcal{H}}(E) = \sup\{s > 0 : \mathcal{H}_s(E) = \infty\} = \inf\{s > 0 : \mathcal{H}_s(E) = 0\}.$$

This is the *Hausdorff dimension* of $E$. Following [10], the *stochastic codimension* of a random set $E$ in $\mathbb{R}^d$, denoted $\text{codim}(E)$, if it exists, is the real number $\beta \in [0, d]$ such that for all compact sets $A \subset \mathbb{R}^d$,

$$\mathbb{P}\{E \cap A \neq \varnothing\} \begin{cases} > 0, & \text{whenever } \dim_{\mathcal{H}}(A) > \beta, \\ = 0, & \text{whenever } \dim_{\mathcal{H}}(A) < \beta. \end{cases}$$

The following result gives a relationship between Hausdorff dimension and stochastic codimension.

THEOREM 3.3 ([10], Theorem 4.7.1, Chapter 11). *Given a random set $E$ in $\mathbb{R}^d$ whose codimension $\beta$ is strictly between $0$ and $d$,*

$$\dim_{\mathcal{H}}(E) + \operatorname{codim}(E) = d \qquad a.s.$$

Theorems 3.1 and 3.3 imply the following result.

COROLLARY 3.4. *Under the hypotheses of Theorem 3.1(b),*

$$\operatorname{codim}\{X((0, +\infty)^N)\} = (d - (N/\alpha))^+$$

*and if $d > N/\alpha$,*

$$\dim_{\mathcal{H}}\{X((0, +\infty)^N)\} = N/\alpha \qquad a.s.$$

PROOF. By Frostman's theorem (see [10], Appendix C, Theorem 2.2.1), the capacitarian and Hausdorff dimensions agree on compact sets. Therefore, the desired result follows from Theorems 3.1 and 3.3.  □

PROOF OF THEOREM 3.1. (a) Suppose $0 < a < b < \infty$ and $M > 0$ are fixed. We assume that Hypothesis A1 holds for $\alpha \in (0,1)$ and $\gamma \geq \alpha$ fixed, and show that Hypotheses H1 and H2 of Theorem 2.4 hold.

VERIFICATION OF HYPOTHESIS H1. Fix $x \in \mathbb{R}^d$ such that $\|x\| \leq M$. Inequality (3.1) implies that

$$
\int_{[a,b]^N} p_{X_t}(x)\, dt = (2\pi)^{-d/2} \int_{[a,b]^N} \sigma^{-d}(t) \exp\left(-\frac{\|x\|^2}{2\sigma^2(t)}\right) dt
$$
$$
\geq (2\pi)^{-d/2} (C_2)^{-d/2} (b-a)^N \exp\left(-\frac{M^2}{2C_1}\right),
$$

which proves Hypothesis H1.

VERIFICATION OF HYPOTHESIS H2. By (3.3) and (3.4), for all $s, t \in [a,b]^N$ with $t \neq s$, $(X_t, X_s)$ has a (Gaussian) density $p_{X_t, X_s}(x, y)$. The latter can be written as

$$p_{X_t, X_s}(x, y) = p_{X_t | X_s = y}(x) p_{X_s}(y),$$

where $p_{X_t | X_s = y}$ denotes the conditional density function of the random variable $X_t$ given $X_s = y$.

Note that the conditional distribution of $X_t^i$ given $X_s^i = y^i$ is Normal with mean $m(s,t)y^i$, where $m(s,t) = (\sigma(s,t))/(\sigma^2(s))$, and variance $\tau^2(s,t) = \sigma^2(t)(1 - \rho^2(s,t))$. Observe that $\tau^2(s,t) > 0$ by (3.3) and (3.4), that (3.2) is a condition on the conditional mean $m(s,t)$ and that (3.3) is a condition on the conditional variance.

Fix $x, y \in \mathbb{R}^d$ such that $\|x\| \leq M$ and $\|y\| \leq M$. By (3.4), $p_{X_t, X_s}(\cdot, \cdot)$ is bounded by some constant $C' > 0$ when $\|t - s\| > \delta$. Therefore,

$$\int_{[a,b]^N} \int_{[a,b]^N} p_{X_t, X_s}(x, y) \, dt \, ds \leq C + \iint_{D(\delta)} p_{X_t, X_s}(x, y) \, dt \, ds,$$

where $D(\delta) = \{(s, t) \in [a, b]^N \times [a, b]^N : \|t - s\| \leq \delta\}$ and $C = C'(b - a)^{2N}$. The integral on the right-hand side can be written

$$(3.7) \qquad (2\pi)^{-d} \iint_{D(\delta)} \tau^{-d}(s, t) \exp\left(-\frac{\|x - m(s,t)y\|^2}{2\tau^2(s,t)}\right) \\ \times \sigma^{-d}(s) \exp\left(-\frac{\|y\|^2}{2\sigma^2(s)}\right) dt \, ds.$$

By the triangle inequality and the identity $(u - v)^2 \geq u^2/2 - v^2$,

$$(3.8) \qquad \exp\left(-\frac{\|x - m(s,t)y\|^2}{2\tau^2(s,t)}\right) \\ \leq \exp\left(-\frac{\|x - y\|^2}{4\tau^2(s,t)}\right) \exp\left(\frac{\|y\|^2|1 - m(s,t)|^2}{2\tau^2(s,t)}\right).$$

By (3.2), there exists a constant $C_3$ such that

$$(3.9) \qquad |1 - m(s,t)|^2 \leq C_3 \|t - s\|^{2\gamma}.$$

By (3.1) and (3.3),

$$(3.10) \qquad K_1 \|t - s\|^{2\alpha} \leq \tau^2(s,t) \leq K_2 \|t - s\|^{2\alpha}.$$

We now apply (3.8)–(3.10) to (3.7). Because $\gamma \geq \alpha$, this yields

$$\int_{[a,b]^N} \int_{[a,b]^N} p_{X_t, X_s}(x, y) \, dt \, ds$$

$$\leq C + (2\pi)^{-d}(K_1)^{-d/2}(C_1)^{-d/2} \iint_{D(\delta)} dt \, ds \|t - s\|^{-\alpha d}$$

$$\times \exp\left(-\frac{\|x - y\|^2}{4K_2\|t - s\|^{2\alpha}}\right) \exp\left(\frac{C_3^2\|y\|^2\|t - s\|^{2\gamma}}{2K_1\|t - s\|^{2\alpha}}\right) \exp\left(-\frac{\|y\|^2}{2C_2}\right)$$

$$\leq C + K_3 \iint_{D(\delta)} \|t - s\|^{-\alpha d} \exp\left(-\frac{\|x - y\|^2}{4K_2\|t - s\|^{2\alpha}}\right) dt \, ds.$$

We now fix $t$ and use the change of variables $u = t - s$ to see that this expression is less than or equal to

$$C + K_3(b - a)^N \int_{B(\delta)} \|u\|^{-\alpha d} \exp\left(-\frac{\|x - y\|^2}{4K_2\|u\|^{2\alpha}}\right) du,$$

where $B(\delta) = \{u \in \mathbb{R}^N : \|u\| \le \delta\}$. Finally, use the change of variables $u = \|x-y\|^{1/\alpha} z(4K_2)^{-1/(2\alpha)}$ to see that this is less than or equal to

$$(3.11) \qquad \begin{aligned} &C + K_4 \|x-y\|^{-d+(N/\alpha)} \\ &\times \int_{B((4K_2)^{1/(2\alpha)}\delta/\|x-y\|^{1/\alpha})} \|z\|^{-\alpha d} \exp(-1/\|z\|^{2\alpha}) \, dz. \end{aligned}$$

We now state a real variable technical lemma which is crucial for our estimates. The proof of this lemma is left to the reader.

LEMMA 3.5. *Define $\varphi_{\alpha,\beta}(r) = \int_{B(r)} \|z\|^{-\beta} e^{-1/\|z\|^{2\alpha}} \, dz$, for all $r > 0$. Then for any $r_0 > 0$ and $\alpha \in (0,1)$, there exist finite constants $c_1, c_2, c_3, c_4 > 0$ such that for all $r \ge r_0$,*

$$c_1 \le \varphi_{\alpha,\beta}(r) \le c_2, \qquad \text{if } \beta > N,$$

$$c_3 \ln(r/r_0) \le \varphi_{\alpha,\beta}(r) \le c_4 \ln(r), \qquad \text{if } \beta = N.$$

Continuing the verification of Hypothesis H2, apply Lemma 3.5 with $\beta = \alpha d$ to (3.11) and use the fact that $C \le C(2M)^{d-(N/\alpha)}\|x-y\|^{-d+(N/\alpha)}$ because $\|x\| \le M$ and $\|y\| \le M$, to conclude the verification of Hypothesis H2 for $d > N/\alpha$. When $d = N/\alpha$, choose $r_0 > 0$ such that $(4K_2)^{1/(2\alpha)}\delta/(2M)^{1/\alpha} \ge r_0$ and apply Lemma 3.5 to (3.11) to obtain

$$\int_{[a,b]^N} \int_{[a,b]^N} p_{X_t, X_s}(x,y) \, dt \, ds \le C + c_4 K_4 \ln\left(\frac{(4K_2)^{1/(2\alpha)}\delta}{\|x-y\|^{1/\alpha}}\right).$$

Note that for all $x, y \in \mathbb{R}^d$ with $\|x\| \le M$ and $\|y\| \le M$, $\|x-y\| \le 2M$. Then we can check that if $\|x-y\| \le 2M$, there exists a finite constant $C' > 1$ such that

$$\ln\left(\frac{(4K_2)^{1/(2\alpha)}\delta}{\|x-y\|^{1/\alpha}}\right) \le C' \ln\left(\frac{3M}{\|x-y\|}\right).$$

On the other hand, note that $C \le C(\ln(3/2))^{-1} \ln(3M/\|x-y\|)$, and the verification of Hypothesis H2 for $d = N/\alpha$ is completed. When $d < N/\alpha$, the expression in (3.11) is bounded, so Hypothesis H2 holds with $k(x) \equiv 1$.

This completes the proof of Theorem 3.1(a).

(b) We now assume that Hypotheses A1 and A2 hold for $\alpha \in (0,1)$ fixed and $\gamma \ge \alpha$, and show that Hypothesis H3 of Theorem 2.4 holds. We also assume that $d \ge N/\alpha$, since otherwise, the statement is trivial.

VERIFICATION OF HYPOTHESIS H3. For all $s, t \in (0, \infty)$ with $s < t$, and for almost all $x \in \mathbb{R}^d$ and $\omega \in \Omega$, let $p_{s,t}(\omega, x)$ denote the conditional density of the random variable $X_t - X_s$ given $\mathcal{F}_s$. By (3.4), the latter exists and, for almost all $\omega \in \Omega$, is a Gaussian density with conditional mean $\mu(s,t) = \mathbb{E}[X_t - X_s | \mathcal{F}_s]$ and deterministic variance $\beta^2(s,t) = \mathbb{E}[(X_t - X_s) - \mathbb{E}[X_t - X_s | \mathcal{F}_s]]^2$ (see [4], Chapter II, Section 3).

It suffices to check the inequality in Hypothesis H3 when $\|x\| \leq 2M$, since the indicator on the right-hand side vanishes for $\|x\| > 2M$. Fix $s \in [a,b]^N$. Then

$$(3.12) \qquad \int_{[b, 2b-a]^N} p_{s,t}(\omega, x) \, dt$$
$$= (2\pi)^{-d/2} \int_{[b,2b-a]^N} \beta^{-d}(s,t) \exp\left(-\frac{\|x - \mu(s,t)\|^2}{2\beta^2(s,t)}\right) dt.$$

By the triangle inequality,

$$(3.13) \quad \exp\left(-\frac{\|x - \mu(s,t)\|^2}{2\beta^2(s,t)}\right) \geq \exp\left(-\frac{\|x\|^2}{2\beta^2(s,t)}\right) \exp\left(-\frac{\|\mu(s,t)\|^2}{2\beta^2(s,t)}\right).$$

We now apply (3.5), (3.6) and (3.13) to (3.12), and we see that this expression is greater than or equal to

$$K_4 \int_{[b,2b-a]^N} \|t-s\|^{-\alpha d} \exp\left(-\frac{\|x\|^2}{C_7 \|t-s\|^{2\alpha}}\right) \exp\left(-\frac{C_6^2 \|t-s\|^{2\gamma}}{C_7 \|t-s\|^{2\alpha}}\right) dt.$$

Using the fact that $\gamma \geq \alpha$ and the change of variables $t-s = \|x\|^{1/\alpha} z (C_7)^{-1/(2\alpha)}$, we get that this is greater than or equal to

$$(3.14) \quad K_5 \|x\|^{-d+(N/\alpha)} \int_{B((C_7)^{1/(2\alpha)}(b-a)/\|x\|^{1/\alpha})} \|z\|^{-\alpha d} \exp(-1/\|z\|^{2\alpha}) \, dz.$$

It now suffices to choose $r_0 > 0$ such that $(C_7)^{1/(2\alpha)}(b-a)/(2M)^{1/\alpha} \geq 3^{1/\alpha} r_0$ and apply Lemma 3.5 to (3.14). This concludes the verification of Hypothesis H3 when $d > N/\alpha$. Finally, for $d = N/\alpha$, we apply Lemma 3.5 to (3.14), to obtain

$$\int_{[b,2b-a]^N} p_{s,t}(\omega, x) \, dt \geq c_3 K_5 \ln\left(\frac{6^{1/\alpha} M^{1/\alpha}}{\|x\|^{1/\alpha}}\right) \geq \frac{1}{\alpha} c_3 K_5 \ln\left(\frac{3M}{\|x\|}\right),$$

which concludes the proof of Hypothesis H3 for $d \geq N/\alpha$.

This completes the proof of Theorem 3.1(b).  □

### 3.2. *Examples.*

3.2.1. *The Brownian sheet.* Suppose that $W = (W_t = (W_t^1, \ldots, W_t^d), t \in \mathbb{R}_+^N)$ is a $d$-dimensional $N$-parameter Brownian sheet, that is, the coordinate processes of $W$ are Gaussian, with zero mean and covariances

$$\mathbb{E}[W_s^k W_t^j] = \prod_{i=1}^{N}(s_i \wedge t_i)\delta_{kj} \qquad \text{for all } s, t \in \mathbb{R}_+^N, 1 \leq k, j \leq d,$$

and defined on the canonical probability space $(\Omega, \mathcal{G}, \mathbb{P})$, where $\Omega$ is the space of all continuous functions $\omega : \mathbb{R}_+^N \to \mathbb{R}^d$ vanishing on the axes, $\mathbb{P}$ is the law of $W$ and $\mathcal{G}$ is the completion of the Borel $\sigma$-field of $\Omega$ with respect to $\mathbb{P}$. We also consider the increasing family of $\sigma$-fields $\mathcal{F} = (\mathcal{F}_t, t \in \mathbb{R}_+^N)$, such that for any $t \in \mathbb{R}_+^N$, $\mathcal{F}_t$ is generated by the random variables $(W_s, s \leq t)$, and the null sets of $\mathcal{G}$. The latter is a complete, right continuous, commuting filtration (see [10], Chapter 7, Theorem 2.4.1).

**THEOREM 3.6** ([11], Theorem 1.1).  *The $(N, d)$ Brownian sheet satisfies Hypotheses* A1 *and* A2 *with $\alpha = 1/2$ and $\gamma = 1$, and therefore the conclusions of Theorem* 3.1.

PROOF.  Fix $0 < a < b < \infty$. Inequality (3.1) is trivial since $\mathbb{E}[W_t^2] = \prod_{i=1}^N t_i$, for all $t \in \mathbb{R}_+^N$. Moreover, since the Brownian sheet has independent increments, (3.5) holds also trivially.

We claim that for all $s, t \in [a, b]^N$ with $s \leq t$, there exist two constants $C_1$ and $C_2$ such that

(3.15)            $$C_1 \|t - s\| \leq \mathbb{E}[(W_t - W_s)^2] \leq C_2 \|t - s\|.$$

Suppose $s \leq t$ and note that $\mathbb{E}[(W_t - W_s)^2] = \sigma^2(t) - \sigma^2(s) = \prod_{i=1}^N t_i - \prod_{i=1}^N s_i$. Let $f(t) = \prod_{i=1}^N t_i$ and let $\phi = (\phi_1, \dots, \phi_N) : [0, 1] \to \mathbb{R}^N$ be an affine function such that $\phi(0) = s$ and $\phi(1) = t$. The partial derivatives of $f$ are bounded above and below on $[a, b]^N$, so

$$f(t) - f(s) = f(\phi(1)) - f(\phi(0)) = \int_0^1 \sum_{i=1}^N \frac{\partial f}{\partial t_i}(\phi(u)) \dot{\phi}_i(u) \, du$$

$$\leq C \sum_{i=1}^N (\phi_i(1) - \phi_i(0)) \leq C' \|t - s\|,$$

which proves the upper bound in (3.15). The lower bound is obtained by proceeding along the same lines.

We now prove (3.2) and (3.3). If $s, t \in [a, b]^N$, using (3.1) and (3.15),

$$\left| \frac{\sigma^2(s) - \sigma(s, t)}{\sigma^2(s)} \right| = \left| \frac{\sigma^2(s) - \sigma^2(s \wedge t)}{\sigma^2(s)} \right| \leq C \|s - (s \wedge t)\| \leq C \|t - s\|,$$

which proves (3.2). To prove (3.3), if $s, t \in [a, b]^N$, using (3.1) and (3.15),

$$\frac{\sigma^2(t)\sigma^2(s) - (\sigma(s, t))^2}{\sigma^2(t)\sigma^2(s)} \geq \frac{\sigma^2(s \wedge t)}{2\sigma^2(t)\sigma^2(s)}(\sigma^2(t) - \sigma^2(s \wedge t) + \sigma^2(s) - \sigma^2(s \wedge t))$$

$$\geq C(\|t - (s \wedge t)\| + \|s - (s \wedge t)\|)$$

$$\geq C \|t - s\|,$$

which concludes the proof of (3.3).

Since (3.6) follows from (3.15), it remains to prove (3.4). Because $s, t \mapsto \rho(s, t)$ is continuous on $[a, b]^N$, it suffices to check that $\rho(s, t) < 1$ for any $s, t \in [a, b]^N$ with $s \neq t$. This is clear since

$$\rho(s, t) = \prod_{i=1}^{N} \frac{s_i \wedge t_i}{\sqrt{s_i t_i}} < 1.$$

This proves Theorem 3.6.  □

3.2.2. $\alpha$-regular Gaussian fields. For $\alpha \in (0, 1)$, an $N$-parameter, $\mathbb{R}^d$-valued process $X = (X_t, t \in \mathbb{R}_+^N)$ is said to be $\alpha$-regular (see [10], Chapter 11, Section 5.2) if $X$ is a centered, stationary Gaussian process with i.i.d. coordinate processes whose covariance function $R$ satisfies the following:

ASSUMPTION R1.    If $s \neq 0$, then $|R(s)| < 1$, and there exist finite positive constants $c_1 \leq c_2$ and $\delta > 0$, such that, for all $s \in \mathbb{R}_+^N$ with $\|s\| \leq \delta$,

$$c_1 \|s\|^{2\alpha} \leq 1 - R(s) \leq c_2 \|s\|^{2\alpha}.$$

Such Gaussian fields were studied in [10], Chapter 11, Section 5.2. In the latter, the stochastic codimension and the Hausdorff dimension of the range of these processes is obtained (see [10], Chapter 11, Theorem 5.2.1 and Corollary 5.2.1).

The following result proves that the $\alpha$-regular Gaussian fields satisfy the lower bound of Theorem 3.1(a). The upper bound cannot be obtained from Theorem 3.1(b) since such processes are not necessarily adapted to a commuting filtration. An upper bound involving Hausdorff measure is given in [10], Chapter 11, Example 5.3.3.

THEOREM 3.7.    Fix $\alpha \in (0, 1)$ with $N/\alpha \geq 2$. Let $X = (X_t, t \in \mathbb{R}_+^N)$ be $\alpha$-regular. Fix $0 < a < b < \infty$ and $0 < M < \infty$. Then there exists a finite positive constant $K$ such that for all compact sets $A \subset \{x \in \mathbb{R}^d : \|x\| < M\}$,

$$K \operatorname{Cap}_{d-(N/\alpha)}(A) \leq \mathbb{P}\{\exists t \in [a, b]^N : X_t \in A\}.$$

PROOF.    Under Assumption R1, $R(0) = 1$, so $\sigma^2(t) = 1$ and $\rho(s, t) = R(t - s) = \sigma(s, t)$. In particular, (3.1) with $\gamma = 2\alpha$ and (3.2) hold. Condition (3.3) holds (for a sufficiently small $\delta > 0$) because

$$1 - \rho^2(s, t) = (1 + \rho(s, t))(1 - \rho(s, t)) \geq c_1(1 - c_2\|t - s\|^{2\alpha})\|t - s\|^{2\alpha}$$
$$\geq c_1(1 - c_2\delta^{2\alpha})\|t - s\|^{2\alpha}.$$

Since $R(\cdot)$ is a bounded covariance function, it is nonnegative definite. By Bochner's theorem ([25], Theorem 6.5.64), $R(\cdot)$ is the Fourier transform of a nonnegative measure, which is a probability measure since $R(0) = 1$. Therefore, $R(\cdot)$ is in fact continuous, so by Assumption R1, there is $\varepsilon > 0$ such that $|R(t-s)| < 1 - \varepsilon$ for $\|t - s\| > \delta$ with $s, t \in [a, b]^N$. Therefore, (3.4) holds. This proves that Hypothesis A1 holds with $\gamma = 2\alpha$, so the conclusion follows from Theorem 3.1(a). $\quad\square$

3.2.3. *The $(N, d)$ Ornstein–Uhlenbeck sheet.* The $(N, d)$ Ornstein–Uhlenbeck sheet $U = (U_t, t \in \mathbb{R}_+^N)$ is defined as

$$U_t = e^{-|t|/2} W_{e^t}, \qquad t \in \mathbb{R}_+^N,$$

where $|t| = |t_1| + \cdots + |t_N|$, $e^t = (e^{t_1}, \ldots, e^{t_N})$ and $W = (W_t, t \in \mathbb{R}_+^N)$ is an $\mathbb{R}^d$-valued Brownian sheet. Therefore, $U$ is an $N$-parameter centered stationary Gaussian process on $\mathbb{R}^d$ with covariance function given by

$$\mathbb{E}[U_s U_t] = e^{-|t-s|/2}, \qquad s, t \in \mathbb{R}_+^N.$$

Consider its natural and completed filtration, denoted $\mathcal{F}_t^U$, that is,

$$\mathcal{F}_t^U = \sigma\{U_s, s \le t\} = \sigma\{W_s, s \le e^t\} = \mathcal{F}_{e^t}^W, \qquad t \in \mathbb{R}_+^N,$$

where $\mathcal{F}_t^W$ denotes the natural filtration of the Brownian sheet and is therefore a commuting filtration.

THEOREM 3.8. *The $(N, d)$ Ornstein–Uhlenbeck sheet satisfies Hypotheses A1 and A2 with $\alpha = 1/2$ and $\gamma = 1$, and therefore the estimates of Theorem 3.1.*

PROOF. Since the $(N, d)$ Ornstein–Uhlenbeck sheet is an $\alpha$-regular Gaussian field with $\alpha = 1/2$, the lower bound follows from Theorem 3.7 and it remains to verify Hypothesis A2. Using the estimate $1 - e^{-x} \le x$ for all $x \ge 0$ and using (3.15), we have, for all $s, t \in [a, b]^N$ with $s \le t$,

$$
\begin{aligned}
(3.16) \quad & |\mathbb{E}[U_t - U_s | \mathcal{F}_s^U]| \\
& = |\mathbb{E}[e^{-|t|/2}(W_{e^t} - W_{e^s}) + (e^{-|t|/2} - e^{-|s|/2})W_{e^s} | \mathcal{F}_{e^s}^W]| \\
& = |(e^{-|t|/2} - e^{-|s|/2})W_{e^s}| \\
& \le |U_s|(|t| - |s|)/2 \\
& \le C_1 M \|t - s\|
\end{aligned}
$$

on $\{|U_s| \le M\}$, which proves (3.5). To prove (3.6), use (3.16) to see that the expectation in (3.6) is equal to

$$\mathbb{E}[e^{-|t|/2} W_{e^t} - e^{-|s|/2} W_{e^s} - (e^{-|t|/2} - e^{-|s|/2})W_{e^s}]^2 = \mathbb{E}[e^{-|t|/2}(W_{e^t} - W_{e^s})]^2,$$

so by (3.15), there exist two constants $C_2$ and $C_3$ such that for all $s, t \in [a, b]^N$ with $s \le t$,

$$C_2 \|e^t - e^s\| \le e^{-|t|} \mathbb{E}[(W_{e^t} - W_{e^s})^2] \le C_3 \|e^t - e^s\|.$$

Finally, using the inequalities $cx \leq 1 - e^{-x} \leq x$, for $x \in [0, b]$, we obtain

$$C_4 \|t - s\| \leq e^{-|t|} \mathbb{E}[(W_{e^t} - W_{e^s})^2] \leq C_5 \|t - s\|,$$

which proves (3.6). The upper bound follows then from Theorem 3.1(b). $\quad\square$

3.2.4. *The fractional Brownian sheet.* The fractional Brownian sheet with Hurst parameter $H \in (0, 1)$, $X^H = (X_t^H, t \in \mathbb{R}_+^N)$, is a $d$-dimensional centered Gaussian process with independent coordinate processes, each with covariance function given by

$$E[(X_t^H)^j (X_s^H)^j] = \prod_{i=1}^N \frac{c}{2}(t_i^{2H} + s_i^{2H} - |t_i - s_i|^{2H}), \qquad s, t \in \mathbb{R}_+^N, 1 \leq j \leq d,$$

where $c$ is a positive finite constant. Note that if $H = \frac{1}{2}$, one obtains the standard Brownian sheet.

As in Theorem 3.7, we obtain only a lower capacity estimate for the fractional Brownian sheet as a consequence of Theorem 3.1(a), since this process is not necessarily adapted to a commuting filtration.

THEOREM 3.9. *Fix $H \in (0, 1)$, $0 < a < b < \infty$ and $0 < M < \infty$. Then there exists a finite positive constant $K$ such that for all compact sets $A \subset \{x \in \mathbb{R}^d : \|x\| < M\}$,*

$$K \operatorname{Cap}_{d-(N/H)}(A) \leq \mathbb{P}\{\exists t \in [a, b]^N : X_t^H \in A\}.$$

PROOF. By Theorem 3.1(a), it suffices to prove that Hypothesis A1 holds with $\alpha = H$ and $\gamma = 2H \wedge 1$. Fix $0 < a < b < \infty$. Since $\sigma^2(t) = c \prod_{i=1}^N t_i^{2H}$ for all $t \in \mathbb{R}_+^N$, (3.1) holds trivially. On the other hand, for all $s, t \in [a, b]^N$, we have

$$(3.17) \quad \left|1 - \frac{\sigma(s, t)}{\sigma^2(s)}\right| = \left|\frac{\prod_{i=1}^N s_i^{2H} - \prod_{i=1}^N ((s_i^{2H} + t_i^{2H} - |t_i - s_i|^{2H})/2)}{\prod_{i=1}^N s_i^{2H}}\right|.$$

To prove (3.2) for $s$ fixed, set $f(u) = \prod_{i=1}^N \frac{1}{2}(s_i^{2H} + u_i^{2H} - |u_i - s_i|^{2H})$ and let $\phi = (\phi_1, \ldots, \phi_N) : [0, 1] \to \mathbb{R}^N$ be an affine function such that $\phi(0) = t$ and $\phi(1) = s$. Because $f$ is differentiable in the orthant centered at $s$ that contains $t$, and in this orthant, for $\|t - s\|$ sufficiently small and $u \in \operatorname{Im} \phi$, $|\frac{\partial f}{\partial u_i}(u)| \leq C(|u_i - s_i|^{2H-1} \vee 1)$, the numerator on the right-hand side of (3.17) can be bounded by

$$|f(\phi(1)) - f(\phi(0))| \leq C \int_0^1 \sum_{i=1}^N (|\phi_i(r) - s_i|^{2H-1} \vee 1) \dot{\phi}_i(r) \, dr = C|t - s|^{2H \wedge 1}.$$

This concludes the proof of (3.2).

To prove (3.3), set $f(u) = \prod_{i=1}^{N} g(u_i)$, where $g(u_i) = \frac{1}{2} u_i^{-H}(1 + u_i^{2H} - |u_i - 1|^{2H})$. An elementary calculation shows that for $|u_i - 1|$ sufficiently small, there are constants $\tilde{c}_2 > \tilde{c}_1 > 0$ such that $\tilde{c}_1 |u_i - 1|^{2H-1} \le |g'(u_i)| \le \tilde{c}_2 |u_i - 1|^{2H-1}$. Therefore, there are constants $c_2 > c_1 > 0$ such that for $\|u - (1, \ldots, 1)\|$ sufficiently small,

$$(3.18) \qquad c_1 |u_i - 1|^{2H-1} \le \left| \frac{\partial f}{\partial u_i}(u) \right| \le c_2 |u_i - 1|^{2H-1}.$$

Let $\phi = (\phi_1, \ldots, \phi_N): [0,1] \to \mathbb{R}^N$ be an affine function such that $\phi(1) = (1, \ldots, 1)$ and $\phi(0) = (t_1/s_1, \ldots, t_N/s_N)$. Observe that

$$1 - \rho(s,t) = f(\phi(1)) - f(\phi(0)) = \int_0^1 \sum_{i=1}^{N} \frac{\partial f}{\partial u_i}(\phi(r)) \dot{\phi}_i(r) \, dr.$$

By (3.18), for $s, t \in [a,b]^2$ with $\|s - t\|$ sufficiently small, this expression is bounded above and below by constants times

$$\sum_{i=1}^{N} \int_0^1 |\phi_i(r) - 1|^{2H-1} \dot{\phi}_i(r) \, dr = \sum_{i=1}^{N} |t_i - s_i|^{2H}/s_i^{2H}.$$

Because $1 - \rho^2(s,t) = (1 - \rho(s,t))(1 + \rho(s,t))$ and the second factor is less than or equal to 2 and bounded away from 0, it follows that (3.3) holds.

It remains to prove (3.4). For this, it suffices to check that for any $s_i, t_i \in [a,b]$ with $t_i - s_i > 0$,

$$t_i^{2H} + s_i^{2H} - |t_i - s_i|^{2H} < 2\sqrt{t_i^{2H} s_i^{2H}}.$$

Move the square root to the left-hand side and isolate the perfect square to see that this is equivalent to $t_i^H - s_i^H < |t_i - s_i|^H$, which holds since $0 < H < 1$. The lower bound follows then from Theorem 3.1(a). □

## 4. Gaussian-type bounds for densities of solutions of hyperbolic SPDEs.
In this section, we use techniques of Malliavin calculus to establish Gaussian-type bounds on the density of the solution to a nonlinear hyperbolic SPDE. These bounds are such that they can be used to verify Hypotheses H1–H3, as we see in Section 5.

4.1. *Elements of Malliavin calculus.* In this section, we recall, following [18], some elements of Malliavin calculus. Let $W = (W_t = (W_t^1, \ldots, W_t^d), t \in \mathbb{R}_+^2)$ be an $\mathbb{R}^d$-valued two-parameter Brownian sheet defined on its canonical probability space $(\Omega, \mathcal{G}, \mathbb{P})$ and let $\mathcal{F} = (\mathcal{F}_t, t \in \mathbb{R}_+^2)$ be its natural filtration (see Section 3.2.1).

Let $H$ be the Hilbert space $H = L^2(\mathbb{R}_+^2, \mathbb{R}^d)$. For any $h \in H$, we set $W(h) = \sum_{j=1}^{d} \int_{\mathbb{R}_+^2} h_j(z) \, dW_z^j$. The Gaussian subspace $\mathcal{H} = \{W(h), h \in H\}$ of $L^2(\Omega, \mathcal{G}, \mathbb{P})$ is isomorphic to $H$.

Let $S$ denote the class of smooth random variables $F = f(W(h_1), \ldots, W(h_n))$, where $h_1, \ldots, h_n$ are in $H$, $n \geq 1$, and $f$ belongs to $C_P^\infty(\mathbb{R}^n)$, the set of functions $f$ such that $f$ and all its partial derivatives have at most polynomial growth.

Given $F$ in $S$, its *derivative* is the $d$-dimensional stochastic process $DF = (D_t F = (D_t^{(1)} F, \ldots, D_t^{(d)} F), t \in \mathbb{R}_+^2)$ given by

$$D_t F = \sum_{i=1}^n \frac{\partial f}{\partial x_i}(W(h_1), \ldots, W(h_n)) h_i(t).$$

More generally, the $k$th-order derivative of $F$ is obtained by iterating the derivative operator $k$ times: if $F$ is a smooth random variable and $k$ is an integer, set $D_{t_1, \ldots, t_k}^k F = D_{t_1} \cdots D_{t_k} F$. Then for every $p \geq 1$ and any natural number $k$, we denote by $\mathbb{D}^{k,p}$ the closure of $S$ with respect to the seminorm $\| \cdot \|_{k,p}$ defined by

$$\|F\|_{k,p}^p = \mathbb{E}[|F|^p] + \sum_{j=1}^k \mathbb{E}[\|D^j F\|_{H^{\otimes j}}^p],$$

where $H^{\otimes j}$ is the product space $L^2((\mathbb{R}_+^2)^j, \mathbb{R}^d)$, and

$$\|D^j F\|_{H^{\otimes j}} = \left( \sum_{k_1, \ldots, k_j = 1}^d \int_{\mathbb{R}_+^2} \overset{(j)}{\cdots} \int_{\mathbb{R}_+^2} |D_{t_1}^{(k_1)} \cdots D_{t_j}^{(k_j)} F|^2 \, dt_1 \cdots dt_j \right)^{1/2}.$$

We set $\mathbb{D}^\infty = \bigcap_{p \geq 1} \bigcap_{k \geq 1} \mathbb{D}^{k,p}$.

Similarly, for any separable Hilbert space $V$, we can define the analogous spaces $\mathbb{D}^{k,p}(V)$ and $\mathbb{D}^\infty(V)$ of $V$-valued random variables, and the related $\| \cdot \|_{k,p,V}$ seminorms (the related smooth functionals being of the form $F = \sum_{j=1}^n F_j v_j$, where $F_j \in S$ and $v_j \in V$).

We denote by $\delta$ the adjoint of the operator $D$, which is an unbounded operator on $L^2(\Omega, H)$ taking values in $L^2(\Omega)$ ([18], Definition 1.3.1). In particular, if $u$ belongs to $\operatorname{Dom} \delta$, then $\delta(u)$ is the element of $L^2(\Omega)$ characterized by the duality relationship

$$(4.1) \qquad \mathbb{E}(F\delta(u)) = \mathbb{E}\left( \sum_{j=1}^d \int_{\mathbb{R}_+^2} D_t^{(j)} F u_t^j \, dt \right) \qquad \text{for any } F \in \mathbb{D}^{1,2}.$$

If $u \in L^2(\mathbb{R}_+^2 \times \Omega, \mathbb{R}^d)$ is an adapted process, then (see [18], Proposition 1.3.4) $u$ belongs to $\operatorname{Dom} \delta$ and $\delta(u)$ coincides with the Itô integral

$$\delta(u) = \sum_{j=1}^d \int_{\mathbb{R}_+^2} u_s^j \, dW_s^j.$$

We use the following estimate of the $\| \cdot \|_{k,p}$ norm of $\delta(u)$.

PROPOSITION 4.1 ([18], Proposition 3.2.1, and [19], (1.11) and page 131). *The adjoint $\delta$ is a continuous operator from $\mathbb{D}^{k+1,p}(H)$ to $\mathbb{D}^{k,p}$ for all $p > 1$, $k \geq 0$. Hence, for all $u \in \mathbb{D}^{k+1,p}(H)$,*

$$\|\delta(u)\|_{k,p} \leq c_{k,p}\|u\|_{k+1,p,H} \tag{4.2}$$

*for some constant $c_{k,p} > 0$.*

Our first application of Malliavin calculus to the study of probability laws is the following global criterion for smoothness of densities.

THEOREM 4.2 ([18], Theorem 2.1.2 and Corollary 2.1.2). *Let $F = (F^1, \dots, F^d)$ be a random vector satisfying the following two conditions:*

(i) *For all $i = 1, \dots, d$, $F^i \in \mathbb{D}^\infty$.*

(ii) *The Malliavin matrix of $F$ defined by $\gamma_F = (\langle DF^i, DF^j \rangle_H)_{1 \leq i,j \leq d}$ is invertible a.s.*

*Then the probability law of $F$ is absolutely continuous with respect to Lebesgue measure. Moreover, assuming* (i) *and* (ii), *and*

(iii) $(\det \gamma_F)^{-1} \in L^p$ *for all $p \geq 1$,*

*the probability density function of $F$ is infinitely differentiable.*

A random vector $F$ that satisfies conditions (i)–(iii) of Theorem 4.2 is said to be *nondegenerate*. For a nondegenerate random vector, the following *integration by parts formula* plays a key role.

PROPOSITION 4.3 ([19], Proposition 3.2.1). *Let $F = (F^1, \dots, F^d) \in (\mathbb{D}^\infty)^d$ be a nondegenerate random vector, let $G \in \mathbb{D}^\infty$ and let $g \in C_p^\infty(\mathbb{R}^d)$. Fix $k \geq 1$. Then for any multiindex $\alpha = (\alpha_1, \dots, \alpha_k) \in \{1, \dots, d\}^k$, there exists an element $H_\alpha(F, G) \in \mathbb{D}^\infty$ such that*

$$\mathbb{E}[(\partial_\alpha g)(F)G] = \mathbb{E}[g(F)H_\alpha(F,G)], \tag{4.3}$$

*where the random variables $H_\alpha(F, G)$ are recursively given by*

$$H_{(i)}(F,G) = \sum_{j=1}^d \delta(G(\gamma_F^{-1})_{ij}DF^j),$$

$$H_\alpha(F,G) = H_{(\alpha_k)}(F, H_{(\alpha_1,\dots,\alpha_{k-1})}(F,G)).$$

4.2. *Conditional Malliavin calculus.* In this section, we give the conditional version of some of the results established in Section 4.1. The proofs are very similar to the one-parameter case and are left to the reader.

The first result is the conditional version of the duality relationship (4.1), which is the two-parameter version of [17], (2.12).

PROPOSITION 4.4. *Let $s, t \in \mathbb{R}_+^2$ with $s < t$. Let $F$ be a random variable in $\mathbb{D}^{1,2}$ and let $u$ be an adapted process such that $\mathbb{E}[\int_{\mathbb{R}_+^2} \|u_s\|^2 ds] < \infty$. Then the following duality relationship holds:*

$$(4.4) \qquad \mathbb{E}\left[F \int_{[0,t] \setminus [0,s]} u_r \, dW_r \Big| \mathcal{F}_s\right] = \mathbb{E}\left[\int_{[0,t] \setminus [0,s]} \langle D_r F, u_r \rangle \, dr \Big| \mathcal{F}_s\right].$$

The following norms are the two-parameter versions of those in [17], Definition 1. Let $s, t \in \mathbb{R}_+^2$ with $s < t$. For any function $f \in L^2([0,t]^n, \mathbb{R}^d)$, any random variable $F \in \mathbb{D}^{k,p}$ and any process $u$ such that $u_r \in \mathbb{D}^{k,p}$ for all $r \in [0,t]$, we define

$$H_{s,t} = L^2([0,t] \setminus [0,s], \mathbb{R}^d),$$

$$\|f\|_{H_{s,t}^{\otimes n}} = \left(\int_{([0,t] \setminus [0,s])^n} \|f(r)\|_{\mathbb{R}^d}^2 \, dr\right)^{1/2},$$

$$\|F\|_{k,p,H_{s,t}}^{\mathcal{F}_s} = \left\{ \mathbb{E}[|F|^p | \mathcal{F}_s] + \sum_{j=1}^k \mathbb{E}[\|D^j F\|_{H_{s,t}^{\otimes j}}^p | \mathcal{F}_s] \right\}^{1/p}$$

and

$$\|u\|_{k,p,H_{s,t}}^{\mathcal{F}_s} = \left\{ \mathbb{E}[\|u\|_{H_{s,t}}^p | \mathcal{F}_s] + \sum_{j=1}^k \mathbb{E}[\|D^j u\|_{H_{s,t}^{\otimes j+1}}^p | \mathcal{F}_s] \right\}^{1/p}.$$

Moreover, we write $\gamma_F^{s,t}$ for the Malliavin covariance matrix with respect to $H_{s,t}$, that is,

$$\gamma_F^{s,t} = (\langle DF^i, DF^j \rangle_{H_{s,t}})_{1 \leq i,j \leq d}.$$

With this notation, we can state the following conditional version of inequality (4.2), which is the two-parameter version of [17], (2.15).

PROPOSITION 4.5. *Let $s, t \in \mathbb{R}_+^2$ with $s < t$. For any $u \in \mathbb{D}^{k+1,p}(H_{s,t})$, we have*

$$(4.5) \qquad \|\delta(u)\|_{k,p,H_{s,t}}^{\mathcal{F}_s} \leq c_{k,p} \|u\|_{k+1,p,H_{s,t}}^{\mathcal{F}_s}$$

*for some constant $c_{k,p} > 0$.*

Finally, we give the conditional version of the two-parameter Burkholder inequality (see [18], A.2, for the two-parameter nonconditional version).

PROPOSITION 4.6. *Fix $p > 1$. There is a finite constant $b_p > 0$ such that for all adapted $X = (X_t, t \in \mathbb{R}_+^2)$ in $L^2(\mathbb{R}_+^2 \times \Omega)$ and all $s, t \in \mathbb{R}_+^2$ with $s < t$,*

$$(4.6) \qquad \mathbb{E}\left[\left|\int_{[0,t] \setminus [0,s]} X_r \, dW_r\right|^p \Big| \mathcal{F}_s\right] \leq b_p \mathbb{E}\left[\left|\int_{[0,t] \setminus [0,s]} X_r^2 \, dr\right|^{p/2} \Big| \mathcal{F}_s\right].$$

4.3. *Hyperbolic stochastic partial differential equations.* Let $b, \sigma_j : \mathbb{R}^d \to \mathbb{R}^d$, $1 \leq j \leq d$, be measurable globally Lipschitz functions, where the vector-valued functions $\sigma_1, \ldots, \sigma_d$ denote the columns of a matrix $\sigma = (\sigma_j^i)_{1 \leq i, j \leq d}$.

Consider the system of stochastic integral equations on the plane

$$(4.7) \quad X_t^i = x_0 + \sum_{j=1}^d \int_{[0,t]} \sigma_j^i(X_s) \, dW_s^j + \int_{[0,t]} b^i(X_s) \, ds, \qquad t \in \mathbb{R}_+^2, 1 \leq i \leq d,$$

where the first integral is an Itô integral with respect to the Brownian sheet (as defined in [28], Chapter 4) and $x_0 \in \mathbb{R}^d$ is the constant value of the process $X_t$ on the axes. It is well known (see [21], Lemma 3.1) that there exists a unique two-parameter, $d$-dimensional, continuous and adapted process $X = (X_t, t \in \mathbb{R}_+^2)$ that satisfies equation (4.7). In addition, $\mathbb{E}[\sup_{r \in [0,t]} |X_r|^p] < \infty$ for any $p \geq 2$ and $t \in \mathbb{R}_+^2$. In [21], Malliavin calculus is used to establish the following result.

THEOREM 4.7 ([21], Proposition 3.3).    *If the coefficients of $\sigma$ and $b$ are infinitely differentiable with bounded partial derivatives of all orders, then $X_t^i$ belongs to $\mathbb{D}^\infty$ for all $t \in \mathbb{R}_+^2$ and $i = 1, \ldots, d$.*

Assuming the latter infinite differentiability condition on the coefficients of $\sigma$ and $b$, we state the following standard *hypoellipticity hypothesis*:

CONDITION P.    There is $n \geq 1$ such that the vector space spanned by the column vectors $\sigma_1, \ldots, \sigma_d$, $\sigma_i \nabla \sigma_j$, $1 \leq i, j \leq d$, $\sigma_i \nabla (\sigma_j \nabla \sigma_k)$, $1 \leq i, j, k \leq d, \ldots, \sigma_{i_1} \nabla (\cdots (\sigma_{i_{n-1}} \nabla \sigma_{i_n}) \cdots)$, $1 \leq i_1, \ldots, i_n \leq d$, at the point $x_0$ is $\mathbb{R}^d$, where the column vector $\sigma_i \nabla \sigma_j$ denotes the covariant derivative of $\sigma_j$ in the direction of $\sigma_i$.

The following result uses Theorem 4.2 and gives the existence and smoothness of the density of $X_t$ for any $t$ away from the axes.

THEOREM 4.8 ([18], Theorem 2.4.2, and [21], Theorem 4.3).    *Under Condition P, for any point $t$ away from the axes, the random vector $X_t$ is nondegenerate and therefore has an absolutely continuous probability distribution with respect to Lebesgue measure on $\mathbb{R}^d$. Moreover, its probability density function is infinitely differentiable.*

4.4. *Gaussian-type upper bounds.* In this section, we present some preliminary results and establish a Gaussian-type upper bound for the drift-free case with vanishing initial conditions.

Let $X = (X_t, t \in \mathbb{R}_+^2)$ be the unique solution of (4.7) with $b \equiv 0$ and $x_0 = 0$, that is,

$$(4.8) \qquad X_t^i = \sum_{j=1}^d \int_{[0,t]} \sigma_j^i(X_s) \, dW_s^j, \qquad t \in \mathbb{R}_+^2, 1 \le i \le d.$$

We assume that the following two hypotheses on the matrix $\sigma$ hold:

HYPOTHESIS P1.  *The coefficients of the matrix $\sigma$ are bounded and infinitely differentiable with bounded partial derivatives (we denote by $T$ the uniform bound on the coefficients of $\sigma$ and its first partial derivatives).*

HYPOTHESIS P2.  Strong ellipticity: $\|\sigma(x)\xi\|^2 = \sum_{k=1}^d (\sum_{i=1}^d \sigma_k^i(x)\xi^i)^2 \ge \rho^2 > 0$ *for some $\rho > 0$, for all $x \in \mathbb{R}^d$ and for all $\xi \in \mathbb{R}^d$ with $\|\xi\| = 1$.*

Note that Hypothesis P2 implies Condition P. Indeed, Hypothesis P2 implies that the vector space spanned by the column vectors $\sigma_1(x), \ldots, \sigma_d(x)$ at any point $x$ in $\mathbb{R}^d$ is $\mathbb{R}^d$, so Condition P holds.

Fix $s \in \mathbb{R}_+^N$. Let $P_s(\omega, \cdot)$ be a regular version of the conditional distribution of the process $(X_t - X_s, t \in \mathbb{R}_+^N \setminus [0, s])$ given $\mathcal{F}_s$, as defined in Section 2. As in Theorem 4.8, we can check that under Condition P for $\mathbb{P}$-almost all $\omega \in \Omega$ and for any $s < t$, the law of $X_t - X_s$ under $P_s(\omega, \cdot)$ is absolutely continuous with respect to Lebesgue measure on $\mathbb{R}^d$. We let $p_{s,t}(\omega, x)$ denote the density of $X_t - X_s$ under $P_s(\omega, \cdot)$. We note that $p_{s,t}(\omega, x)$ is the conditional density of $X_t - X_s$ given $\mathcal{F}_s$. Therefore, for a random variable $Y$, we interpret $\mathbb{E}[Y|\mathcal{F}_s]$ as $E_s(\omega, Y)$, where $E_s(\omega, \cdot)$ denotes the expectation under $P_s(\omega, \cdot)$. However, for $s \in \mathbb{R}_+^N$ fixed, there is $N_s \in \mathcal{F}_s$ with $\mathbb{P}(N_s) = 0$ such that $p_{s,t}(\omega, \cdot)$ is defined for $\omega \in \Omega \setminus N_s$ and all $s < t$.

LEMMA 4.9.  *Assume Hypotheses* P1 *and* P2. *Fix* $s, t \in \mathbb{R}_+^2$ *with $s < t$, and $x \in \mathbb{R}^d$. Let $\sigma$ be a subset of the set of indices $\{1, \ldots, d\}$. Then, for all $\omega \in \Omega \setminus N_s$ and every $\sigma$,*

$$p_{s,t}(\omega, x) = (-1)^{d-|\sigma|} \mathbb{E}[\mathbb{1}_{\{X_t^i - X_s^i > x^i, \, i \in \sigma, \, X_t^i - X_s^i < x^i, \, i \notin \sigma\}} H_{(1,\ldots,d)}(X_t - X_s, 1)|\mathcal{F}_s],$$

*where $|\sigma|$ is the cardinality of $\sigma$ and, in agreement with Proposition* 4.3,

$$H_{(1,\ldots,d)}(X_t - X_s, 1)$$
$$= \delta(((\gamma_{X_t - X_s}^{s,t})^{-1} D(X_t - X_s))^d \delta(\cdots \delta((\gamma_{X_t - X_s}^{-1} D(X_t - X_s))^1) \cdots)).$$

Lemma 4.9 is a consequence of the integration by parts formula (4.3). It gives an explicit expression for the density of the process that will be very useful for further computations. When $\sigma = \{1, \ldots, d\}$, the proof of Lemma 4.9 follows along the same lines as the one-parameter nonconditional case (see [19], Corollary 3.2.1) and is therefore omitted. For the general case, see [22], (5.3), where a similar expression is obtained for the kernel of a stochastic semigroup.

A key property is that the moments of the iterated derivatives of $X$ are finite; see [17], Lemma 6, where similar nonconditional estimates are obtained for one-parameter Brownian martingales in $\mathbb{R}^d$.

LEMMA 4.10. *Assuming Hypothesis* P1, *for any* $0 < a < b < \infty$, $p \geq 1$ *and* $n \geq 1$, *there exists a finite constant* $C > 0$ *depending on* $a, b$ *and the uniform bounds from Hypothesis* P1 *such that for any* $s, t \in [a, b]^2$ *with* $s < t$, *on* $\Omega \setminus N_s$,

$$(4.9) \qquad \sup_{\substack{z_1 \leq r, \dots, z_n \leq r \\ z_1, \dots, z_n, r \in [0, t] \setminus [0, s]}} \mathbb{E}[|D_{z_1}^{(k_1)} \cdots D_{z_n}^{(k_n)}(X_r^i)|^p | \mathcal{F}_s] \leq C$$

*for* $1 \leq i, k_1, \dots, k_n \leq d$.

PROOF. We prove this lemma by induction on $n$. Suppose $n = 1$. For any $z \in [0, t] \setminus [0, s]$ and $s, t \in \mathbb{R}_+^2$ with $s < t$, the process $(D_z^{(k)}(X_t^i), 1 \leq k \leq d)$ satisfies the following system of stochastic differential equations (see [18], page 127):

$$(4.10) \qquad D_z^{(k)}(X_t^i) = \sigma_k^i(X_z) + \sum_{j=1}^d \int_{[z, t]} D_z^{(k)} \sigma_j^i(X_r) \, dW_r^j.$$

Using Burkholder's inequality (4.6) for conditional expectations, for any $p \geq 1$,

$$\sum_{i=1}^d \mathbb{E}[|D_z^{(k)}(X_t^i)|^p | \mathcal{F}_s]$$

$$\leq \sum_{i=1}^d 2^{p-1} \left\{ \mathbb{E}[|\sigma_k^i(X_z)|^p | \mathcal{F}_s] + \mathbb{E}\left[ \left| \sum_{j=1}^d \int_{[z, t]} D_z^{(k)} \sigma_j^i(X_r) \, dW_r^j \right|^p \Big| \mathcal{F}_s \right] \right\}$$

$$\leq 2^{p-1} \left\{ dT^p + b_p d^{p-1} \sum_{i, j=1}^d \mathbb{E}\left[ \left| \int_{[z, t]} |D_z^{(k)} \sigma_j^i(X_r)|^2 \, dr \right|^{p/2} \Big| \mathcal{F}_s \right] \right\}$$

$$\leq 2^{p-1} \left\{ dT^p + b_p d^{2p} b^{p-2} T^p \int_{[z, t]} \sum_{i=1}^d \mathbb{E}[|D_z^{(k)}(X_r^i)|^p | \mathcal{F}_s] \, dr \right\}.$$

Finally, for $z$ fixed, using a two-parameter version of Gronwall's lemma (see [18], Exercise 2.4.3) in the form

$$(4.11) \qquad f(z, t) \leq A + B \int_{[z, t]} f(z, r) \, dr,$$

we conclude the proof of (4.9) for $n = 1$.

We now assume that (4.9) holds for $n > 1$. Apply $n$ times the derivative operator to (4.10), which yields

$$D_{z_1}^{(k_1)} \cdots D_{z_{n+1}}^{(k_{n+1})}(X_t^i) = \sum_{l=1}^{n+1} D_{z_1}^{(k_1)} \cdots D_{z_{l-1}}^{(k_{l-1})} D_{z_{l+1}}^{(k_{l+1})} \cdots D_{z_{n+1}}^{(k_{n+1})}(\sigma_{k_l}^i(X_{z_l}))$$

$$+ \sum_{j=1}^{d} \int_{[z_1 \vee \cdots \vee z_{n+1}, t]} D_{z_1}^{(k_1)} \cdots D_{z_{n+1}}^{(k_{n+1})}(\sigma_j^i(X_r)) \, dW_r^j.$$

For the first term, we use the induction hypothesis and the uniform bounds on the derivatives of the coefficients of $\sigma$. For the second term, we use Burkholder's inequality and again the induction hypothesis and the bounds on the derivatives on the coefficients of $\sigma$. We finally obtain

$$\sum_{i=1}^{d} \mathbb{E}[|D_{z_1}^{(k_1)} \cdots D_{z_{n+1}}^{(k_{n+1})}(X_t^i)|^p]$$

$$\leq C_1 + C_2 \mathbb{E}\left[\int_{[z_1 \vee \cdots \vee z_{n+1}, t]} \sum_{i=1}^{d} |D_{z_1}^{(k_1)} \cdots D_{z_{n+1}}^{(k_{n+1})}(X_r^i)|^p \, dr\right]$$

and the proof is completed using Gronwall's lemma.  □

In the next lemma, we follow [17], Lemma 12, where similar estimates are carried out for one-parameter Brownian martingales in $\mathbb{R}^d$.

LEMMA 4.11. *Assuming Hypotheses* P1 *and* P2, *there exists a finite constant* $C > 0$ *depending on* $a, b$ *and the uniform bounds from Hypotheses* P1 *and* P2 *such that, for any* $s, t \in [a, b]^2$ *with* $s < t$, *on* $\Omega \setminus N_s$,

$$(4.12) \qquad (\mathbb{E}[\{H_{(1,\ldots,d)}(X_t - X_s, 1)\}^2 | \mathcal{F}_s])^{1/2} \leq C\|t - s\|^{-d/2}.$$

PROOF. To simplify the notation, we write $H_{(1,\ldots,d)}^{t,s} = H_{(1,\ldots,d)}(X_t - X_s, 1)$. Using (4.5) and Hölder's inequality for conditional Sobolev norms in Wiener space (see [29], Proposition 1.10, page 50), we obtain

$$\|H_{(1,\ldots,d)}^{t,s}\|_{0,2,H_{s,t}}^{\mathcal{F}_s}$$

$$= \|H_{(d)}(X_t - X_s, H_{(1,\ldots,d-1)}^{t,s})\|_{0,2,H_{s,t}}^{\mathcal{F}_s}$$

$$= \left\|\sum_{j=1}^{d} \delta(H_{(1,\ldots,d-1)}^{t,s}((\gamma_{X_t-X_s}^{s,t})^{-1})_{dj} D(X_t^j - X_s^j))\right\|_{0,2,H_{s,t}}^{\mathcal{F}_s}$$

$$\leq c\|H_{(1,\ldots,d-1)}^{t,s}\|_{1,4,H_{s,t}}^{\mathcal{F}_s} \sum_{j=1}^{d} \|((\gamma_{X_t-X_s}^{s,t})^{-1})_{dj}\|_{1,8,H_{s,t}}^{\mathcal{F}_s} \|D(X_t^j - X_s^j)\|_{1,8,H_{s,t}}^{\mathcal{F}_s}.$$

Hence, to prove (4.12), it suffices to show that for each $p \geq 1$ and $n = 1, \ldots, d$,

$$(4.13) \qquad \|D(X_t^j - X_s^j)\|_{n,p,H_{s,t}}^{\mathcal{F}_s} \leq c_{n,p}^1 \|t - s\|^{1/2}$$

for some finite constant $c_{n,p}^1 > 0$ and

$$(4.14) \qquad \|((\gamma_{X_t-X_s}^{s,t})^{-1})_{ij}\|_{n,p,H_{s,t}}^{\mathcal{F}_s} \leq c_{n,p}^2 \|t - s\|^{-1}$$

for some finite constant $c_{n,p}^2 > 0$.

Indeed, (4.13) and (4.14) with $n = 1$ and $p = 8$ imply that

$$\|H_{(1,\ldots,d)}^{t,s}\|_{0,2,H_{s,t}}^{\mathcal{F}_s} \leq c d c^1 c^2 \|t - s\|^{-1/2} \|H_{(1,\ldots,d-1)}^{t,s}\|_{1,4,H_{s,t}}^{\mathcal{F}_s}.$$

Iterating the process, we find

$$\|H_{(1,\ldots,d)}^{t,s}\|_{0,2,H_{s,t}}^{\mathcal{F}_s} \leq (c d c^1 c^2)^d \|t - s\|^{-d/2},$$

which concludes the proof of (4.12). □

PROOF OF (4.13).  Fix $p \geq 1$. By definition,

$$(4.15) \qquad \begin{aligned} &\|D(X_t^i - X_s^i)\|_{n,p,H_{s,t}}^{\mathcal{F}_s} \\ &= \left\{ \mathbb{E}[\|D(X_t^i - X_s^i)\|_{H_{s,t}}^p | \mathcal{F}_s] + \sum_{k=2}^{n+1} \mathbb{E}[\|D^k(X_t^i - X_s^i)\|_{H_{s,t}^{\otimes k}}^p | \mathcal{F}_s] \right\}^{1/p}. \end{aligned}$$

Furthermore, for any $z \in [0,t] \setminus [0,s]$, the process $(D_z^{(k)}(X_t^i - X_s^i), 1 \leq k \leq d)$ satisfies the system of stochastic differential equations

$$(4.16) \qquad D_z^{(k)}(X_t^i - X_s^i) = \sigma_i^k(X_z) + \sum_{j=1}^d \int_{[z,t]} D_z^{(k)} \sigma_j^i(X_r) \, dW_r^j.$$

By Burkholder's inequality (4.6) for conditional expectations and using the Cauchy–Schwarz inequality, for any $p \geq 1$,

$$\begin{aligned} &\mathbb{E}[\|D(X_t^i - X_s^i)\|_{H_{s,t}}^p | \mathcal{F}_s] \\ &\quad = \mathbb{E}\left[ \left| \sum_{k=1}^d \int_{[0,t] \setminus [0,s]} |D_z^{(k)}(X_t^i - X_s^i)|^2 \, dz \right|^{p/2} \Big| \mathcal{F}_s \right] \\ &\quad \leq d^{p/2-1} (c\|t - s\|)^{p/2-1} \sum_{k=1}^d \mathbb{E}\left[ \int_{[0,t] \setminus [0,s]} |D_z^{(k)}(X_t^i - X_s^i)|^p \, dz \, \Big| \mathcal{F}_s \right] \\ &\quad \leq 2^{p-1} d^{p/2-1} (c\|t - s\|)^{p/2-1} \\ &\qquad \times \left\{ d T^p c \|t - s\| \right. \end{aligned}$$

$$+ d^{p-1} \sum_{k,j=1}^{d} \mathbb{E}\left[\int_{[0,t]\setminus[0,s]}\left|\int_{[z,t]} D_z^{(k)}\sigma_j^i(X_r)\,dW_r^j\right|^p dz\,\Big|\mathcal{F}_s\right]\Bigg\}$$

$$\leq 2^{p-1}(c\|t-s\|)^{p/2}$$

$$\times \Bigg\{ d^{p/2}T^p + (c\|t-s\|)^{p/2}d^{5p/2-2}b_p(T)^p$$

$$\times \sum_{k,l=1}^{d} \sup_{r,z\in[0,t]\setminus[0,s]} \mathbb{E}[|D_z^{(k)}(X_r^l)|^p|\mathcal{F}_s]\Bigg\}.$$

Finally, by Lemma 4.10, we get that there exists a positive finite constant $k_1(a,b,T)$ such that

$$(4.17) \qquad \mathbb{E}[\|D(X_t^i - X_s^i)\|_{H_{s,t}}^p|\mathcal{F}_s] \leq k_1\|t-s\|^{p/2}.$$

To estimate the second term on the right-hand side of (4.15), we apply $n$ times the derivative operator to (4.16) to get

$$D_{z_1}^{(k_1)}\cdots D_{z_{n+1}}^{(k_{n+1})}(X_t^i - X_s^i) = \sum_{l=1}^{n+1} D_{z_1}^{(k_1)}\cdots D_{z_{l-1}}^{(k_{l-1})}D_{z_{l+1}}^{(k_{l+1})}\cdots D_{z_{n+1}}^{(k_{n+1})}(\sigma_{k_l}^i(X_{z_l}))$$

$$+ \sum_{j=1}^{d}\int_{[z_1\vee\cdots\vee z_{n+1},t]} D_{z_1}^{(k_1)}\cdots D_{z_{n+1}}^{(k_{n+1})}(\sigma_j^i(X_r))\,dW_r^j.$$

Proceeding as above, by Burkholder's inequality (4.6) for conditional expectations, and using the Cauchy–Schwarz inequality and Lemma 4.10, we obtain, for all $k = 2, \ldots, n+1$,

$$(4.18) \qquad \mathbb{E}[\|D^k(X_t^i - X_s^i)\|_{H_{s,t}^{\otimes k}}^p|\mathcal{F}_s] \leq k_2\|t-s\|^{p/2},$$

where $k_2$ is a finite positive constant.

Finally, substituting (4.17) and (4.18) into (4.15) concludes the proof of (4.13). □

PROOF OF (4.14). Fix $p \geq 1$. By definition

$$(4.19) \qquad \begin{aligned} &\|((\gamma_{X_t-X_s}^{s,t})^{-1})_{ij}\|_{n,p,H_{s,t}}^{\mathcal{F}_s} \\ &= \Bigg\{ \mathbb{E}[(|(\gamma_{X_t-X_s}^{s,t})^{-1})_{ij}|^p|\mathcal{F}_s] \\ &\quad + \sum_{k=1}^{n} \mathbb{E}[\|D^k((\gamma_{X_t-X_s}^{s,t})^{-1})_{ij}\|_{H_{s,t}^{\otimes k}}^p|\mathcal{F}_s]\Bigg\}^{1/p}. \end{aligned}$$

We use a standard argument to estimate the moments of the inverse of the Malliavin matrix. We follow [18], proof of (3.22), and [17], Lemma 10. Using the Cauchy–Schwarz inequality and Cramér's formula for $(\gamma_{X_t-X_s}^{s,t})^{-1}$, we can easily check that for all $p \geq 1$,

$$
\begin{aligned}
(4.20) \quad & \mathbb{E}[(((\gamma_{X_t-X_s}^{s,t})^{-1})_{ij})^p|\mathcal{F}_s] \\
& \leq c_{d,p}\mathbb{E}[(\det \gamma_{X_t-X_s}^{s,t})^{-2p}|\mathcal{F}_s]^{1/2} \times \mathbb{E}[\|D(X_t-X_s)\|_{H_{s,t}}^{4p(d-1)}|\mathcal{F}_s]^{1/2}
\end{aligned}
$$

for some constant $c_{d,p} > 0$. For the second factor, we use (4.17) to get

$$
(4.21) \qquad \mathbb{E}[\|D(X_t-X_s)\|_{H_{s,t}}^{4p(d-1)}|\mathcal{F}_s] \leq k_1\|t-s\|^{2p(d-1)}
$$

for some finite constant $k_1(a,b,T) > 0$. On the other hand, we write

$$
\begin{aligned}
\det \gamma_{X_t-X_s}^{s,t} &\geq \inf_{\|v\|=1}(v^T\gamma_{X_t-X_s}^{s,t}v)^d \\
&= \inf_{\|v\|=1}\left(\sum_{k=1}^d\int_{[0,t]\setminus[0,s]}\left|\sum_{i=1}^d D_z^{(k)}(X_t^i-X_s^i)v_i\right|^2 dz\right)^d.
\end{aligned}
$$

Using (4.16) and Hypothesis P2, for any $h \in (0,1]$, we see that the expression in parentheses is bounded below by

$$
\begin{aligned}
&\sum_{k=1}^d\int_{[0,t-(1-h)(t-s)]\setminus[0,s]} dz \left|\sum_{i=1}^d v_i\left(\sigma_k^i(X_z)+\sum_{j=1}^d\int_{[z,t]} D_z^{(k)}\sigma_j^i(X_r)\,dW_r^j\right)\right|^2 \\
&\geq \tfrac{1}{2}A\rho^2 - I_h,
\end{aligned}
$$

where $A$ denotes the area of the region $[0,t-(1-h)(t-s)]\setminus[0,s]$ and

$$
I_h = \sum_{k=1}^d\int_{[0,t-(1-h)(t-s)]\setminus[0,s]}\left|\sum_{i,j=1}^d\int_{[z,t]} D_z^{(k)}\sigma_j^i(X_r)\,dW_r^j\right|^2 dz.
$$

We choose $y$ such that $A\rho^2 = 4y^{-1/d}$ and notice that since $h \leq 1$, $y \geq c := 4^d(b\sqrt{2})^{-d}\|t-s\|^{-d}\rho^{-2d}$. In addition, as $h$ varies in $(0,1]$, $y$ varies in $[c,\infty)$. Applying Chebyshev's inequality for conditional probabilities, we find that for any $q \geq 2$,

$$
\begin{aligned}
\mathbb{P}\left\{\det\gamma_{X_t-X_s}^{s,t} < \frac{1}{y}\Big|\mathcal{F}_s\right\} &\leq \mathbb{P}\left\{\left(\frac{1}{2}A\rho^2-I_h\right)^d < \frac{1}{y}\Big|\mathcal{F}_s\right\} \\
&\leq \mathbb{P}\{I_h > y^{-1/d}|\mathcal{F}_s\} \leq y^{q/d}\mathbb{E}[|I_h|^q|\mathcal{F}_s].
\end{aligned}
$$

Using Burkholder's inequality (4.6) for conditional expectations, for any $q \geq 2$, we have

$$
\mathbb{E}[|I_h|^q|\mathcal{F}_s] \leq d^{5q-3}b_qA^{2(q-1)}\sum_{i,j,k=1}^d\mathbb{E}\left[\iint_{([0,t]\setminus[0,s])^2}|D_z^{(k)}\sigma_j^i(X_r)|^{2q}\,dr\,dz\Big|\mathcal{F}_s\right].
$$

By Lemma 4.10(i), the conditional expectation of the right-hand side is bounded above by some finite positive constant $k_2(a, b, T)$. Using the definition of $A$, we obtain

$$\mathbb{E}[|I_h|^q | \mathcal{F}_s] \leq k_3 \frac{4^{2(q-1)}}{\rho^{4(q-1)}} y^{2(1-q)/d}.$$

Consequently, taking $q > 2 + 2pd$,

$$\mathbb{E}[(\det \gamma_{X_t - X_s}^{s,t})^{-2p} | \mathcal{F}_s]$$

$$= \int_0^\infty 2py^{2p-1} \mathbb{P}\{(\det \gamma_{X_t - X_s}^{s,t})^{-1} > y | \mathcal{F}_s\} \, dy$$

$$\leq c^{2p} + 2p \int_c^\infty y^{2p-1} \mathbb{P}\left\{\det \gamma_{X_t - X_s}^{s,t} < \frac{1}{y} \Big| \mathcal{F}_s\right\} dy$$

$$\leq \frac{4^{2dp}}{\|t-s\|^{2dp}(b\sqrt{2})^{2dp}\rho^{4dp}} + 2p \int_c^\infty y^{2p-1+(q/d)} \mathbb{E}[|I_h|^q | \mathcal{F}_s] \, dy$$

$$\leq \frac{4^{2dp}}{\|t-s\|^{2dp}(b\sqrt{2})^{2dp}\rho^{4dp}} + 2pk_3 \frac{4^{2(q-1)}}{\rho^{4(q-1)}} \int_c^\infty y^{2p-1-(q/d)+(2/d)} \, dy$$

$$\leq k_4 \|t-s\|^{-2dp},$$

where $k_4$ is a finite positive constant. Therefore, we have proved that

$$(4.22) \qquad \mathbb{E}[(\det \gamma_{X_t - X_s}^{s,t})^{-2p} | \mathcal{F}_s] \leq k_4 \|t-s\|^{-2dp}.$$

Substituting (4.21) and (4.22) into (4.20), we obtain

$$(4.23) \qquad \begin{aligned} \mathbb{E}[|((\gamma_{X_t - X_s}^{s,t})^{-1})_{ij}|^p | \mathcal{F}_s] &\leq k_5 (\|t-s\|^{-2pd} \|t-s\|^{2p(d-1)})^{1/2} \\ &= k_5 \|t-s\|^{-p} \end{aligned}$$

for some finite constant $k_5 > 0$ not depending on $i$ or $j$. This proves the desired estimate for the first term in (4.19). Turning to the second term, we claim that for all $i, j = 1, \ldots, d$ and $k = 1, \ldots, n$,

$$(4.24) \qquad \mathbb{E}[\|D^k((\gamma_{X_t - X_s}^{s,t})^{-1})_{ij}\|_{H_{s,t}^{\otimes k}}^p | \mathcal{F}_s] \leq k_6 \|t-s\|^{-p}$$

for some finite constant $k_6 > 0$. Indeed, by iterating the equality (see [18], Lemma 2.1.6)

$$D(\gamma_{X_t}^{-1})_{ij} = -\sum_{k,l=1}^d (\gamma_{X_t}^{-1})_{ik} D(\gamma_{X_t})_{kl} (\gamma_{X_t}^{-1})_{jl},$$

and using Hölder's inequality for conditional expectations, we have

$$\sup_{i,j} \mathbb{E}[\|D^k((\gamma_{X_t - X_s}^{s,t})^{-1})_{ij}\|_{H_{s,t}^{\otimes k}}^p | \mathcal{F}_s]$$

$$\leq c \sup_{r=1}^{k} \sum_{\substack{k_1 + \cdots + k_r = k \\ k_l \geq 1, l = 1, \ldots, r}} \mathbb{E}[\|D^{k_1}(\gamma^{s,t}_{X_t - X_s})_{i_1 j_1}\|^{p(r+1)}_{H^{\otimes k_1}_{s,t}} | \mathcal{F}_s]^{1/(r+1)} \times \cdots$$

$$\times \mathbb{E}[\|D^{k_r}(\gamma^{s,t}_{X_t - X_s})_{i_r j_r}\|^{p(1+r)}_{H^{\otimes k_r}_{s,t}} | \mathcal{F}_s]^{1/(r+1)}$$

$$\times \sup_{i,j} \mathbb{E}[|((\gamma^{s,t}_{X_t - X_s})^{-1})_{ij}|^{p(r+1)^2} | \mathcal{F}_s]^{1/(r+1)},$$

where the supremum before the summation is over $i_1, j_1, \ldots, i_r, j_r \in \{1, \ldots, d\}$.

By (4.23), for all $i, j = 1, \ldots, d$,

$$\mathbb{E}[|((\gamma^{s,t}_{X_t - X_s})^{-1})_{ij}|^{p(r+1)^2} | \mathcal{F}_s] \leq k_7 \|t - s\|^{-p(r+1)^2}$$

for some finite constant $k_7 > 0$.

For the other factors, express $D^k(\gamma^{s,t}_{X_t - X_s})_{ij}$ using the definition of $\gamma^{s,t}_{X_t - X_s}$ and use the Cauchy–Schwarz inequality twice and (4.18) to get

$$\mathbb{E}[\|D^k(\gamma^{s,t}_{X_t - X_s})_{ij}\|^p_{H^{\otimes k}_{s,t}} | \mathcal{F}_s]$$

$$= \mathbb{E}\left[\left\|D^k\left(\int_{[0,t]\setminus[0,s]} D_r(X^i_t - X^i_s) \cdot D_r(X^j_t - X^j_s) \, dr\right)\right\|^p_{H^{\otimes k}_{s,t}} \Big| \mathcal{F}_s\right]$$

$$\leq (k+1)^{p-1} \sum_{l=0}^{k} \binom{k}{l}^p \mathbb{E}\left[\left\|\int_{[0,t]\setminus[0,s]} D^l D_r(X^i_t - X^i_s)\right.\right.$$

$$\left.\left. \cdot D^{k-l} D_r(X^j_t - X^j_s) \, dr\right\|^p_{H^{\otimes k}_{s,t}} \Big| \mathcal{F}_s\right]$$

$$\leq d^{p/2}(k+1)^{p-1} \sum_{l=0}^{k} \binom{k}{l}^p \{(\mathbb{E}[\|D^l D(X^i_t - X^i_s)\|^{2p}_{H^{\otimes(l+1)}_{s,t}} | \mathcal{F}_s])^{1/2}$$

$$\times (\mathbb{E}[\|D^{k-l} D(X^j_t - X^j_s)\|^{2p}_{H^{\otimes(k-l+1)}_{s,t}} | \mathcal{F}_s])^{1/2}\}$$

$$\leq k_8 \|t - s\|^p$$

for some finite constant $k_8 > 0$ not depending on $i$ and $j$. This concludes the proof of (4.24). Finally, substituting (4.23) and (4.24) into (4.19), we conclude the proof of (4.14).

The proof of Lemma 4.11 is now complete.    □

The following result is a consequence of Lemmas 4.9 and 4.11.

LEMMA 4.12.    *Assuming Hypotheses* P1 *and* P2, *for any* $0 < a < b < \infty$, *the density function* $p_{X_t}(x)$ *is uniformly bounded for* $t \in [a, b]^2$ *and* $x \in \mathbb{R}^d$.

PROOF. Using Lemma 4.9 with $s = 0$ and $\sigma = \{1, \ldots, d\}$ and the Cauchy–Schwarz inequality, we obtain

$$p_{X_t}(x) \leq (\mathbb{E}[\{H_{(1,\ldots,d)}(X_t, 1)\}^2])^{1/2}.$$

By Lemma 4.11 with $s = 0$, there exists a finite positive constant $c_2$ depending on $a$, $b$ and the uniform bounds from Hypotheses P1 and P2 such that

$$\mathbb{E}[\{H_{(1,\ldots,d)}(X_t, 1)\}^2] \leq c_2,$$

which proves the lemma. $\square$

The next proposition is the main result of this section.

PROPOSITION 4.13. *Assuming Hypotheses* P1 *and* P2, *for any* $0 < a < b < \infty$, *there exists a finite positive constant* $c$ *depending on* $a, b$ *and the uniform bounds from Hypotheses* P1 *and* P2 *such that for any* $x \in \mathbb{R}^d$, $s, t \in [a, b]^2$ *with* $s < t$ *and for* $\omega \in \Omega \setminus N_s$,

$$(4.25) \qquad p_{s,t}(\omega, x) \leq c\|t - s\|^{-d/2} \exp\left(-\frac{\|x\|^2}{c\|t - s\|}\right).$$

PROOF. Apply the Cauchy–Schwarz and Hölder inequalities for conditional expectations to the expression of Lemma 4.9 with $\sigma = \{i \in \{1, \ldots, d\} : x^i \geq 0\}$ to find that

$$(4.26) \qquad p_{s,t}(\omega, x) \leq \prod_{i=1}^{d} (\mathbb{P}\{|X_t^i - X_s^i| \geq |x^i| \,|\, \mathcal{F}_s\})^{1/(2d)} (\mathbb{E}[\{H_{(1,\ldots,d)}^{t,s}\}^2])^{1/2}.$$

Consider the one-parameter martingale $[M_u = (M_u^1, \ldots, M_u^d), 0 \leq u \leq |t| - |s|]$ defined by

$$M_u^i = \begin{cases} \displaystyle\sum_{j=1}^{d} \int_{[s_1, s_1+u] \times [0, s_2]} \sigma_j^i(X_r) \, dW_r^j, & \text{if } 0 \leq u \leq t_1 - s_1, \\[3mm] \displaystyle M_{t_1-s_1}^i + \sum_{j=1}^{d} \int_{[0,t_1] \times [s_2, u+s_2+s_1-t_1]} \sigma_j^i(X_r) \, dW_r^j, \\[3mm] \hspace{4cm} \text{if } t_1 - s_1 \leq u \leq |t| - |s|, \end{cases}$$

for all $i = 1, \ldots, d$, with respect to the filtration $(\mathcal{G}_u, 0 \leq u \leq |t| - |s|)$ defined by

$$\mathcal{G}_u = \begin{cases} \mathcal{F}_{(s_1+u, s_2)}, & \text{if } 0 \leq u \leq t_1 - s_1, \\[2mm] \mathcal{F}_{(t_1, u+s_2+s_1-t_1)}, & \text{if } t_1 - s_1 \leq u \leq |t| - |s|. \end{cases}$$

Notice that $M_0 = 0$, $M_{|t|-|s|} = X_t - X_s$ and $\mathcal{G}_0 = \mathcal{F}_s$. By [1], (2.9),

$$\langle M^i \rangle_{|t|-|s|} = \langle M^i \rangle_{t_1-s_1} + (\langle M^i \rangle_{|t|-|s|} - \langle M^i \rangle_{t_1-s_1})$$

$$= \sum_{j=1}^{d} \int_{[0,t] \setminus [0,s]} \sigma_j^i(X_r)^2 \, dr.$$

Moreover, Hypothesis P1 and the Cauchy–Schwarz inequality imply that $\langle M^i \rangle_{|t|-|s|} \leq C \|t - s\|$, where $C = bd2^{1/2}T^2$ for all $1 \leq i \leq d$. Applying the exponential martingale inequality [18], A.2, we get

$$(4.27) \quad \mathbb{P}\{|X_t^i - X_s^i| \geq |x^i| | \mathcal{F}_s\} \leq 2 \exp\left(-\frac{|x^i|^2}{2C \|t - s\|}\right), \qquad 1 \leq i \leq d.$$

Finally, substituting (4.27) into (4.26) and using Lemma 4.11 concludes the proof of (4.25). $\quad\square$

4.5. *Gaussian-type lower bounds.* In this section we present a lower bound of Gaussian type for the density of the random variable $X_t$ for any $t$ away from the axes, where $X = (X_t, t \in \mathbb{R}_+^2)$ denotes the solution of (4.8), and we present an analogous lower bound for the conditional density of $X_t - X_s$ given $\mathcal{F}_s$, when $s < t$. These results are an application of results by Kohatsu-Higa [13], where Gaussian-type lower bounds are obtained for the density of a general class of uniformly elliptic random variables on a Wiener space, which generalize the lower bound estimates for uniformly elliptic diffusion processes obtained by Kusuoka and Stroock [14].

THEOREM 4.14. *Assume Hypotheses* P1 *and* P2 *and let* $X = (X_t, t \in \mathbb{R}_+^2)$ *be the solution of* (4.8). *Then for any* $0 < a < b < \infty$, *there is a constant* $C$ *which depends only on* $a, b$ *and the uniform bounds from Hypotheses* P1 *and* P2 *such that, for all* $s = (s_1, s_2) \in [a, b]^2$,

$$p_{X_s}(x) \geq C(s_1 s_2)^{-d/2} \exp\left(-\frac{\|x\|^2}{Cs_1 s_2}\right).$$

PROOF. Following [13], Theorem 5, we need to prove that for fixed $s = (s_1, s_2) \in [a, b]^2$, the $d$-dimensional random vector $X_s$ is uniformly elliptic, with constants that do not depend on $s$. Following the same notation as in the Main setup of [13], we replace their $[0, T] \times A$ with $[0, s_1] \times [0, s_2]$ and we set $g(\cdot, \cdot) \equiv 1$. The underlying one-parameter filtration is defined as $(\mathcal{F}_t^1, 0 \leq t \leq s_1)$, where

$$\mathcal{F}_t^1 = \sigma\{W(z_1, z_2), (z_1, z_2) \in [0, t] \times [0, s_2]\}.$$

(In this proof, $t$ denotes a real number, as in [13].)

Consider a sufficiently fine partition $\{0 = t_0 < \cdots < t_N = s_1\}$ of the interval $[0, s_1]$. To simplify the notation, we write $R_{t_{n-1}, t_n}$ for the rectangle $[t_{n-1}, t_n] \times [0, s_2]$. For any $\theta = (\theta_1, \theta_2) \in R_{t_{n-1}, t_n}$, we write $\pi\theta$ for the point $(t_{n-1}, \theta_2)$, that is, the orthogonal projection of $\theta$ on the vertical line through $(t_{n-1}, 0)$. We write $R_{\pi\theta, \theta}$ for the rectangle $[t_{n-1}, \theta_1] \times [0, \theta_2]$. Finally, for any $\theta \in R_{t_{n-1}, t_n}$ and $\eta \in R_{\pi\theta, \theta}$, we write $\theta \triangle \eta$ for the point, that is, the orthogonal projection of $\eta$ on the horizontal line through $\pi\theta$ and $\theta$. Let $\Delta_{n-1}(g)$ be the quantities denoted $\Delta_{i-1}(g)$ in the Main setup of [13], that is,

$$\Delta_{n-1}(g) = \int_{t_{n-1}}^{t_n} \|g(t, \cdot)\|_{L^2([0,s_2],\mathbb{R}^d)}^2 \, dt = (t_n - t_{n-1}) s_2$$

for $n = 1, \ldots, N$. Note that $\|g\|_{L^2([0,s_1] \times [0,s_2], \mathbb{R}^d)} = s_1 s_2$.

We now define the sequence of nondegenerate random vectors $F_n$ required in the Main setup of [13] by $F_n = X_{t_n, s_2}$, for $0 \leq n \leq N$, with $F_N = X_{s_1, s_2}$. Notice that $F_n = (F_n^1, \ldots, F_n^d)$ and $F_n^i = X_{t_n, s_2}^i$. By (4.8),

$$(4.28) \qquad F_n^i - F_{n-1}^i = \sum_{j=1}^d \int_{R_{t_{n-1}, t_n}} \sigma_j^i(X_\theta) \, dW_\theta^j, \qquad 1 \leq i \leq d.$$

The first objective is to find the Itô expansion of order $l \geq 1$ of the random variable $F_n^i - F_{n-1}^i$ to obtain the approximations $\bar{F}_n^i$ required in [13], Theorem 5. For this, we need to introduce some notation. We define the multiindex $\beta \in \bigcup_{n \geq 1} \{0, 1\}^n$, with length $l(\beta)$, and write $\beta = (\beta_1, \ldots, \beta_{l(\beta)})$. We write $-\beta$ and $\overline{\beta}-$ for the multiindex obtained by deleting the first and the last component, respectively, of the multiindex $\beta$. For completeness, we write $\{v\}$ for the multiindex of length zero. These multiindices will be used to write multiple integrals, some of which are stochastic and some are deterministic: when the index is 1, $d^1 W_\theta^j$ denotes $dW_\theta^j$, and when the index is 0, $d^0 W_\theta^j$ denotes $d\theta$. Usually, these integrals are not taken on the whole space $[0, s_1] \times [0, s_2]$, but on subsets of it. The subset most often is of the type $[t_{n-1}, t_n] \times [0, s_2]$. In general, these integrals are denoted by $I_\beta(h_\beta)$ for an $\mathcal{F}_{t_{n-1}}^1$-measurable random process $h_\beta$ such that

$$\sup_{\omega \in \Omega} \|h_\beta\|_{L^2(([t_{n-1}, t_n] \times [0, s_2])^{l(\beta)}, \mathbb{R}^d)}(\omega) < \infty.$$

Given a family of $M$ functions $f_m : \mathbb{R}^d \to \mathbb{R}$, $1 \leq m \leq M$, and $M$ points $\theta^1, \ldots, \theta^M$ placed on the same vertical line in $R_{t_{n-1}, t_n}$, we define for any $1 \leq r \leq d$ and $\xi \in R_{t_{n-1}, t_n}$ the following operations on random variables:

$$L_{r, \theta^1, \ldots, \theta^M, \xi}^1 \left( \prod_{m=1}^M f_m(X_{\theta^m}) \right)$$

$$= \sum_{m=1}^{M} \mathbb{1}_{R_{\pi\theta^m, \theta^m}}(\xi) \left( \prod_{n=1, n\neq m}^{M} f_n(X_{\theta^n \triangle \xi}) \right) \sum_{k=1}^{d} \frac{\partial f_m}{\partial x_k}(X_{\theta^m \triangle \xi}) \sigma_r^k(X_\xi),$$

$$L_{r,\theta^1,\dots,\theta^M,\xi}^0 \left( \prod_{m=1}^{M} f_m(X_{\theta^m}) \right)$$

$$= \left\{ \sum_{m=1}^{M} \frac{1}{2} \mathbb{1}_{R_{\pi\theta^m, \theta^m}}(\xi) \left( \prod_{n=1, n\neq m}^{M} f_n(X_{\theta^n \triangle \xi}) \right) \right.$$

$$\left. \times \sum_{k,l=1}^{d} \frac{\partial^2 f_m}{\partial x_k \, \partial x_l}(X_{\theta^m \triangle \xi}) \sigma_r^k(X_\xi) \sigma_r^l(X_\xi) \right\}$$

$$+ \left\{ \sum_{m_1, m_2=1}^{M} \frac{1}{2} \mathbb{1}_{R_{\pi\theta^{m_1}, \theta^{m_1}} \cap R_{\pi\theta^{m_2}, \theta^{m_2}}}(\xi) \left( \prod_{n=1, n\neq m_1, m_2}^{M} f_n(X_{\theta^n \triangle \xi}) \right) \right.$$

$$\left. \times \sum_{k,l=1}^{d} \frac{\partial f_{m_1}}{\partial x_k}(X_{\theta^{m_1} \triangle \xi}) \sigma_r^k(X_\xi) \frac{\partial f_{m_2}}{\partial x_l}(X_{\theta^{m_2} \triangle \xi}) \sigma_r^l(X_\xi) \right\}.$$

Note that the points $\theta^1 \triangle \xi, \dots, \theta^M \triangle \xi$ and $\xi$ are also placed on a single vertical line in $R_{t_{n-1}, t_n}$.

These operations will be used to apply the standard multidimensional Itô formula to $f(X_{\theta^1}, \dots, X_{\theta^M})$, when $f(x^1, \dots, x^M)$ is of the form $\prod_{m=1}^{M} f_m(x^m)$, where the $f_m$ are as above, $1 \leq m \leq M$. Indeed,

$$u_1 \mapsto (X_{u_1, \theta_2^1}, \dots, X_{u_1, \theta_2^M}), \qquad t_{n-1} \leq u_1 \leq t_n,$$

is an $Md$-dimensional martingale with mutual covariation

$$d\langle X_{\cdot, \theta_2^n}^k, X_{\cdot, \theta_2^m}^l \rangle_{u_1} = du_1 \sum_{r=1}^{d} \int_0^{\theta_2^n \wedge \theta_2^m} \sigma_r^k(X_{u_1, \xi_2}) \sigma_r^l(X_{u_1, \xi_2}) \, d\xi_2.$$

Therefore, according to the Itô formula [2],

$$\prod_{m=1}^{M} f_m(X_{\theta^m}) = \prod_{m=1}^{M} f_m(X_{\pi\theta^m})$$

(4.29)
$$+ \sum_{r=1}^{d} \int_{R_{t_{n-1}, t_n}} L_{r,\theta^1,\dots,\theta^M,\xi}^1 \left( \prod_{m=1}^{M} f_m(X_{\theta^m}) \right) dW_\xi^r$$

$$+ \sum_{r=1}^{d} \int_{R_{t_{n-1}, t_n}} L_{r,\theta^1,\dots,\theta^M,\xi}^0 \left( \prod_{m=1}^{M} f_m(X_{\theta^m}) \right) d\xi.$$

We are now going to iteratively apply the Itô formula to the integrands in (4.29), to obtain an Itô expansion similar to the one presented in [12], Chapter 5. For this, we introduce some additional notation. Given a multiindex

$\beta$, we define the functions

$$h^i_{\beta,j}(\theta) = \sigma^i_j(X_\theta), \qquad\qquad \text{if } \beta = \{v\},$$

$$h^i_{\beta,j,r}(\theta,\xi) = L^1_{r,\theta,\xi}(\sigma^i_j(X_\theta)), \qquad \text{if } \beta = (1),$$

$$h^i_{\beta,j,r}(\theta,\xi) = L^0_{r,\theta,\xi}(\sigma^i_j(X_\theta)), \qquad \text{if } \beta = (0),$$

$$h^i_{\beta,j,r_1,\ldots,r_{l(\beta)}}(\theta,\xi^1,\ldots,\xi^{l(\beta)}) = L^\beta_{r_1,\ldots,r_{l(\beta)},\theta,\xi^1,\ldots,\xi^{l(\beta)}}(\sigma^i_j(X_\theta)),$$

where

$$L^\beta_{r_1,\ldots,r_{l(\beta)},\theta,\xi^1,\ldots,\xi^{l(\beta)}} = L^{\beta_1}_{r_{l(\beta)},\theta,\xi^1,\ldots,\xi^{l(\beta)}} \circ L^{\beta_2}_{r_{l(\beta)-1},\theta,\xi^1,\ldots,\xi^{l(\beta)-1}} \circ \cdots \circ L^{\beta_{l(\beta)}}_{r_1,\theta,\xi^1}.$$

We define the multiple (Itô stochastic/deterministic) integrals $I_\beta$ by

$$
\begin{aligned}
(4.30) \qquad & I_\beta[h^i_{\beta-j,r_1,\ldots,r_{l(\beta)-1}}(\theta,\xi^1,\ldots,\xi^{l(\beta)-1})]_{R_{t_{n-1},t_n}} \\
& = \sum_{j,r_1,\ldots,r_{l(\beta)-1}=1}^d \int_{R_{t_{n-1},t_n}} \Bigl(\int_{R_{\pi\theta,\theta}} \Bigl(\cdots \Bigl(\int \\
& \qquad\qquad \times h^i_{\beta-j,r_1,\ldots,r_{l(\beta)-1}}(\theta,\xi^1,\ldots,\xi^{l(\beta)-1}) \\
& \qquad\qquad\qquad \times d^{\beta_1} W^{r_{l(\beta)-1}}_{\xi^{l(\beta)-1}}\Bigr)\cdots\Bigr) d^{\beta_{l(\beta)-1}} W^{r_1}_{\xi^1}\Bigr) dW^j_\theta,
\end{aligned}
$$

where the domain of integration of each integral is given by the indicator functions that appear in the definition of the operators $L^0$ and $L^1$.

Note that the functions $h^i_{\beta-j,r_1,\ldots,r_{l(\beta)-1}}(\theta,\xi^1,\ldots,\xi^{l(\beta)-1})$ are sums of functions of the form $\prod_{n=1}^{l(\beta)} f_n(X_{\theta^n})$, where $f_n:\mathbb{R}^d \to \mathbb{R}$ are derivatives of some order of the coefficients of the matrix $\sigma$, and $\theta^1,\ldots,\theta^{l(\beta)}$ are $l(\beta)$ points placed on a single vertical line in $R_{t_{n-1},t_n}$ that depends on $\theta,\xi^1,\ldots,\xi^{l(\beta)-1}$. We now define multiple Itô stochastic integrals $I^\pi_\beta$, by modifying the integrand in the definition of $I_\beta$: each point in (4.30) is projected onto the vertical line through $(t_{n-1},0)$, that is,

$$
\begin{aligned}
& I^\pi_\beta[h^i_{\beta-j,r_1,\ldots,r_{l(\beta)-1}}(\theta,\xi^1,\ldots,\xi^{l(\beta)-1})]_{R_{t_{n-1},t_n}} \\
& = \sum_{j,r_1,\ldots,r_{l(\beta)-1}=1}^d \int_{R_{t_{n-1},t_n}} \Bigl(\int_{R_{\pi\theta,\theta}} \Bigl(\cdots\Bigl(\int \sum \prod_{n=1}^{l(\beta)} f_n(X_{\pi\theta^n}) \\
& \qquad\qquad\qquad\qquad \times d^{\beta_1} W^{r_{l(\beta)-1}}_{\xi^{l(\beta)-1}}\Bigr)\cdots\Bigr) \\
& \qquad\qquad\qquad\qquad\qquad \times d^{\beta_{l(\beta)-1}} W^{r_1}_{\xi^1}\Bigr) dW^j_\theta.
\end{aligned}
$$

LEMMA 4.15. *The Itô expansion of order $l \geq 1$,*

$$
\begin{aligned}
F_n^i - F_{n-1}^i &= \sum_{\beta \in \mathcal{A}_l} I_\beta^\pi [h_{\beta-j,r_1,\ldots,r_{l(\beta)-1}}^i (\theta, \xi^1, \ldots, \xi^{l(\beta)-1})]_{R_{t_{n-1},t_n}} \\
&\quad + \sum_{\beta \in \mathcal{B}_l} I_\beta [h_{\beta-j,r_1,\ldots,r_{l(\beta)-1}}^i (\theta, \xi^1, \ldots, \xi^{l(\beta)-1})]_{R_{t_{n-1},t_n}},
\end{aligned}
$$
(4.31)

*holds, where the sets $\mathcal{A}_l$ and $\mathcal{B}_l$ are recursively defined by $\mathcal{A}_1 = \{(1)\}$, $\mathcal{B}_1 = \{(1,1),(0,1)\}$, $\mathcal{A}_{l+1} = \{\mathcal{A}_l, \mathcal{B}_l\}$ and $\mathcal{B}_{l+1} = \{\beta : -\beta \in \mathcal{B}_l\}$ for $l \geq 1$.*

PROOF. We prove this lemma by induction on $l$. By the standard multi-dimensional Itô formula, applied to the random variable $\sigma_j^i(X_\theta)$ with respect to the first coordinate, it follows from (4.28) that

$$
\begin{aligned}
F_n^i - F_{n-1}^i &= \sum_{j=1}^d \int_{R_{t_{n-1},t_n}} \sigma_j^i(X_\theta)\, dW_\theta^j \\
&= \sum_{j=1}^d \int_{R_{t_{n-1},t_n}} \sigma_j^i(X_{\pi\theta})\, dW_\theta^j \\
&\quad + \sum_{j,k,r=1}^d \int_{R_{t_{n-1},t_n}} \left( \int_{R_{\pi\theta,\theta}} \frac{\partial \sigma_j^i(X_{\theta \triangle \eta})}{\partial x_k} \sigma_r^k(X_\eta)\, dW_\eta^r \right) dW_\theta^j \\
&\quad + \frac{1}{2} \sum_{j,k,l,r=1}^d \int_{R_{t_{n-1},t_n}} \left( \int_{R_{\pi\theta,\theta}} \frac{\partial^2 \sigma_j^i(X_{\theta \triangle \eta})}{\partial x_k \partial x_l} \sigma_r^k(X_\eta) \sigma_r^l(X_\eta)\, d\eta \right) dW_\theta^j
\end{aligned}
$$

for all $1 \leq i \leq d$, which proves (4.31) for $l = 1$.

Now assume that (4.31) holds with $\mathcal{A}_{l-1}$ and $\mathcal{B}_{l-1}$ for $l > 1$. Then we apply the Itô formula (4.29) to all the random variables of the form $\prod_{n=1}^{l(\beta)} f_n(X_{\theta^n})$ that appear in the stochastic integrals of the second term of (4.31) with $\mathcal{B}_{l-1}$. Here, $M = l(\beta)$, the functions $f_m : \mathbb{R}^d \to \mathbb{R}$, $1 \leq m \leq l(\beta)$, are derivatives of order greater than or equal to 0 of the coefficients of the matrix $\sigma$ and $\theta^1, \ldots, \theta^{l(\beta)}$ are $l(\beta)$ points placed on a single vertical line in $R_{t_{n-1},t_n}$. The first term on the right-hand side of (4.29) gives a new $I_\beta^\pi$ integral and the two other terms each give a new $I_\beta$ integral. Note that integrals with the operator $L^1$ add a 1 to all the multiindexes of $\mathcal{B}_{l-1}$ and that integrals with the operator $L^0$ add a 0. This yields (4.31) for $l$ with $\mathcal{A}_l = \{\mathcal{A}_{l-1}, \mathcal{B}_{l-1}\}$ and $\mathcal{B}_l = \{\beta : -\beta \in \mathcal{B}_{l-1}\}$. $\quad \square$

Continuing the proof of Theorem 4.14, with this result we define, following [13], Theorem 5, the approximation $\bar{F}_n^i$ of order $l \geq 1$ by

$$
\begin{aligned}
\bar{F}_n^i &= ((t_n - t_{n-1})s_2)^{(l+1)/2} Z_n^i + F_{n-1}^i \\
&\quad + \sum_{\beta \in \mathcal{A}_l} I_\beta^\pi [h_{\beta-j,r_1,\ldots,r_{l(\beta)-1}}^i (\theta, \xi^1, \ldots, \xi^{l(\beta)-1})]_{R_{t_{n-1},t_n}},
\end{aligned}
$$

where $(Z_n, n = 1, \ldots, N)$ is an i.i.d. sequence of $d$-dimensional $N(0, I)$ random variables independent of the Wiener process $W$. In the setting of [13], Theorem 5, we have $k = 1$ and set

$$G_n^l = \sum_{\beta \in \mathcal{A}_l \setminus \mathcal{A}_1} I_\beta^\pi [h_{\beta-j, r_1, \ldots, r_{l(\beta)-1}}^i (\theta, \xi^1, \ldots, \xi^{l(\beta)-1})]_{R_{t_{n-1}, t_n}}.$$

To simplify notation we write $h_\beta^i$ for

$$h_{\beta-j, r_1, \ldots, r_{l(\beta)-1}}^i (\theta, \xi^1, \ldots, \xi^{l(\beta)-1}).$$

With these definitions and since by Hypothesis P1, for any multiindex $\beta \in \mathcal{A}_l$,

$$\sup_{\omega \in \Omega} \|h_\beta^i\|_{L^2(([t_{n-1}, t_n] \times [0, s_2])^{l(\beta)}, \mathbb{R}^d)}(\omega) < \infty,$$

(H1) of [13], Theorem 5 is satisfied. To prove (H2a) of [13], Theorem 5, note that by definition

$$\|F_n^i - \bar{F}_n^i\|_{n, p}^{\mathcal{F}_{t_{n-1}}^1} = \left\| ((t_n - t_{n-1})s_2)^{(l+1)/2} Z_n^i + \sum_{\beta \in \mathcal{B}_l} I_\beta [h_\beta^i]_{R_{t_{n-1}, t_n}} \right\|_{n, p}^{\mathcal{F}_{t_{n-1}}^1}$$

and, for each $\beta \in \mathcal{B}_l$,

$$
\begin{aligned}
(4.32) \quad & \|I_\beta [h_\beta^i]_{R_{t_{n-1}, t_n}}\|_{q, p, H_{t_{n-1}, t_n}}^{\mathcal{F}_{t_{n-1}}^1} \\
& = \left\{ \mathbb{E}[|I_\beta [h_\beta^i]_{R_{t_{n-1}, t_n}}|^p | \mathcal{F}_{t_{n-1}}^1] \right. \\
& \quad \left. + \sum_{j=1}^q \mathbb{E}\left[ \|D^j I_\beta [h_\beta^i]_{R_{t_{n-1}, t_n}}\|_{H_{t_{n-1}, t_n}^{\otimes j}}^p | \mathcal{F}_{t_{n-1}}^1 \right] \right\}^{1/p},
\end{aligned}
$$

where $H_{t_{n-1}, t_n} = L^2([t_{n-1}, t_n] \times [0, s_2], \mathbb{R}^d)$.

Fix $\beta \in \mathcal{B}_l$ and recall that $l(\beta) = l+1$. To estimate the first term in (4.32), we use Burkholder's inequality (4.6) for conditional expectations and the uniform bounds on the derivatives of the coefficients of $\sigma$ to get

$$(4.33) \quad \mathbb{E}[|I_\beta [h_\beta^i]_{R_{t_{n-1}, t_n}}|^p | \mathcal{F}_{t_{n-1}}^1] \leq C((t_n - t_{n-1})s_2)^{(l+1)p/2}$$

for some positive finite constant $C > 0$ independent of the partition, $t_1$, $s$ and $\omega \in \Omega$.

For the second term in (4.32) fix $j \in \{1, \ldots, q\}$. Using formula (1.46) from [18], we can easily check that the term $D^j I_\beta [h_\beta^i]_{R_{t_{n-1}, t_n}}$ contains multiple (Itô stochastic/deterministic) integrals of orders $l+1, \ldots, l+1-j$. Therefore, again using Burkholder's inequality for conditional expectations and the uniform bounds on the derivatives of the coefficients of $\sigma$, we obtain

$$(4.34) \quad \mathbb{E}\left[ \|D^j I_\beta [h_\beta^i]_{R_{t_{n-1}, t_n}}\|_{H_{t_{n-1}, t_n}^{\otimes j}}^p | \mathcal{F}_{t_{n-1}}^1 \right] \leq C((t_n - t_{n-1})s_2)^{(l+1)p/2}$$

for some positive finite constant $C > 0$ independent of the partition, $t_1$, $s$ and $\omega \in \Omega$.

Finally, using (4.33) and (4.34), we get that there exists a constant $C$ independent of the partition, $t_1$, $s$ and $\omega \in \Omega$, such that

$$\|F_n^i - \bar{F}_n^i\|_{n,p}^{\mathcal{F}_{t_{n-1}}^1} \leq C((t_n - t_{n-1})s_2)^{(l+1)/2},$$

which proves (H2a) of [13], Theorem 5, with $\gamma = 1/2$.

Now, using exactly the same argument that led to (4.22) we can easily check that

$$(\mathbb{E}[(\det \gamma_{F_n}^{t_{n-1},t_n})^{-p}|\mathcal{F}_{t_{n-1}}^1])^{1/p} \leq C((t_n - t_{n-1})s_2)^{-d},$$

where $\gamma_{F_n}^{t_{n-1},t_n} = (\langle DF_n^i, DF_n^j \rangle_{L^2([t_{n-1},t_n] \times [0,s_2],\mathbb{R}^d)})_{1 \leq i,j \leq d}$, for some positive finite constant $C$ independent of the partition, $t_1$ and $\omega \in \Omega$. Therefore, (H2b) of [13], Theorem 5, is proved.

Condition (H2c) of [13], Theorem 5, can be rewritten in our case as

$$C_2 \leq ((t_n - t_{n-1})s_2)^{-1} \int_{R_{t_{n-1},t_n}} \sum_{k=1}^d \left( \sum_{i=1}^d \sigma_k^i(X_{\pi\theta})\xi^i \right)^2 d\theta \leq C_1$$

for $\xi \in \mathbb{R}^d$ with $\|\xi\| = 1$. This is obviously satisfied by Hypotheses P1 and P2.

Finally, condition (H2d) of [13], Theorem 5, is satisfied using the same argument that led to condition (H2a), because the higher-order integrals are, in the $\mathcal{F}_{t_{n-1}}^1$ conditional $\mathbb{D}^{n,p}$ norm, smaller than $(t_n - t_{n-1})s_2$. Theorem 4.14 is now proved. □

The next result is the conditional version of Theorem 4.14. Let $p_{s,t}(\omega, x)$ be the density of $X_t - X_s$ under the conditional distribution $P_s(\omega, \cdot)$ defined just above Lemma 4.9, for $\omega \in \Omega \setminus N_s$.

THEOREM 4.16. *Assuming Hypotheses* P1 *and* P2*, for any* $0 < a < b < \infty$*, there exists a positive finite constant* $c$*, which depends on* $a, b$ *and the uniform bounds from Hypotheses* P1 *and* P2*, with the following property: For all* $s, t \in [a, b]^2$ *with* $s < t$*,* $x \in \mathbb{R}^d$*, and for* $\omega \in \Omega \setminus N_s$*,*

$$p_{s,t}(\omega, x) \geq c\|t - s\|^{-d/2} \exp\left(-\frac{\|x\|^2}{c\|t - s\|}\right).$$

PROOF. Fix $s, t \in [a, b]^2$ with $s < t$. As in Theorem 4.14, Theorem 4.16 is an application of [13], Theorem 5. We replace their $[0, T]$ by the union of two segments with extremities $(s_1, s_2)$, $(t_1, s_2)$ and $(t_1, s_2)$, $(t_1, t_2)$, and we set $g(\cdot, \cdot) \equiv 1$. We consider a sufficiently fine partition $\{s = t_0 < \cdots < t_N = t\}$ of the union of the two segments, where $t_i = (t_i^1, t_i^2)$, $1 \le i \le N$. Note that there exist two constants $c_1(a, b)$ and $c_2(a, b)$ such that

$$c_1^{-1}\|t_n - t_{n-1}\| \le \Delta_{n-1}(g) \le c_1\|t_n - t_{n-1}\|$$

and

$$c_2^{-1}\|t - s\| \le \|g\|_{L^2([0,t]\setminus[0,s],\mathbb{R}^d)} \le c_2\|t - s\|.$$

In this case, all of the Wiener stochastic integrals are taken on subsets of one of the $[0, t_n] \setminus [0, t_{n-1}]$. Then, the projection of a given point $\theta \in [0, t_n] \setminus [0, t_{n-1}]$ depends on its position. Namely, if $\theta$ is in the rectangle $[s_1, t_1] \times [0, s_2]$, $\pi\theta$ will be its orthogonal projection on the vertical line through $(t_{n-1}^1, 0)$, while if $\theta$ is in the rectangle $[0, t_1] \times [s_2, t_2]$, $\pi\theta$ will be its orthogonal projection on the horizontal line through $(0, t_{n-1}^2)$.

Define $F_n^i = X_{t_n}^i - X_s^i$ for $1 \le i \le d$, $0 \le n \le N$. Then

$$F_n^i - F_{n-1}^i = \sum_{j=1}^d \int_{[0,t_n]\setminus[0,t_{n-1}]} \sigma_j^i(X_\theta) \, dW_\theta^j, \qquad 1 \le i \le d.$$

Now, using the same notation as in the proof of Theorem 4.14 and proceeding exactly along the same lines, we can easily establish the following result:

LEMMA 4.17. *The Itô expansion of order $l \ge 1$,*

$$F_n^i - F_{n-1}^i = \sum_{\beta \in \mathcal{A}_l} I_\beta[h_{\beta-,j,r_1,\dots,r_{l(\beta)-1}}^i(\theta, \xi_1, \dots, \xi_{l(\beta)-1})]_{[0,t_n]\setminus[0,t_{n-1}]}$$

$$+ \sum_{\beta \in \mathcal{B}_l} I_\beta^\pi[h_{\beta-,j,r_1,\dots,r_{l(\beta)-1}}^i(\theta, \xi_1, \dots, \xi_{l(\beta)-1})]_{[0,t_n]\setminus[0,t_{n-1}]},$$

*holds, where the sets $\mathcal{A}_l$ and $\mathcal{B}_l$ are recursively defined by $\mathcal{A}_1 = \{(1)\}$, $\mathcal{B}_1 = \{(1,1),(0,1)\}$, $\mathcal{A}_{l+1} = \{\mathcal{A}_l, \mathcal{B}_l\}$ and $\mathcal{B}_{l+1} = \{\beta : -\beta \in \mathcal{B}_l\}$ for $l \ge 1$.*

Again, following [13], Theorem 5, we define the approximation $\bar{F}_n^i$ of order $l \ge 1$ by

$$\bar{F}_n^i = \|t_n - t_{n-1}\|^{(l+1)/2} Z_n^i + F_{n-1}^i$$

$$+ \sum_{\beta \in \mathcal{A}_l} I_\beta^\pi[h_{\beta-,j,r_1,\dots,r_{l(\beta)-1}}^i(\theta, \xi_1, \dots, \xi_{l(\beta)-1})]_{[0,t_n]\setminus[0,t_{n-1}]},$$

where $(Z_n^i, n = 0, \dots, N)$ is an i.i.d. sequence of $d$-dimensional $N(0, I)$ random variables independent of the Wiener process $W$.

The remainder of the proof follows exactly as in the proof of Theorem 4.14. We note that the constant $c$ in the conclusion of Theorem 4.16 is also uniform in $\omega \in \Omega \setminus N_s$, $s, t \in [a, b]^2$ with $s < t$ and $x \in \mathbb{R}^d$.    $\square$

**5. Potential theory for hyperbolic SPDEs.**    In this section, we extend the results obtained in Section 2 to the solution of equation (4.7). The proofs make use of Malliavin calculus and are an application of the results of Section 4, which contains the technical work.

5.1. *The case $b \equiv 0$.*    Let $X = (X_t, t \in \mathbb{R}_+^2)$ be the unique $d$-dimensional adapted continuous process defined on $(\Omega, \mathcal{G}, \mathbb{P})$ that solves (4.8). The aim of this section is to establish the following result.

THEOREM 5.1.    *Assuming Hypotheses* P1 *and* P2, *for all* $0 < a < b < \infty$ *and* $M > 0$, *there exists a positive finite constant $K$ depending on $a, b, M, \rho$ and the uniform bounds on the coefficients of $\sigma$ and its derivatives, such that for all compact sets $A \subset \{x \in \mathbb{R}^d : \|x\| < M\}$,*

$$K^{-1} \operatorname{Cap}_{d-4}(A) \leq \mathbb{P}\{\exists t \in [a, b]^2 : X_t \in A\} \leq K \operatorname{Cap}_{d-4}(A),$$

*where* $\operatorname{Cap}_\beta$ *denotes the capacity with respect to the Newtonian $\beta$ kernel $k_\beta(\cdot)$ defined in Theorem 3.1.*

PROOF.    To prove Theorem 5.1 it suffices to prove that Hypotheses H1–H3 of Theorem 2.4 hold for the process $X$. Since Hypothesis H3 is only used for the upper bound in the statement of Theorem 5.1, we only prove Hypothesis H3 when $d \geq 4$, since the upper bound is trivially satisfied for $d \leq 3$. Fix $0 < a < b < \infty$, $M > 0$ and recall that Hypotheses P1 and P2 are satisfied.

VERIFICATION OF HYPOTHESIS H1.    Fix $x \in \mathbb{R}^d$ such that $\|x\| \leq M$. Using Theorem 4.14, we find that

$$\int_{[a,b]^2} p_{X_t}(x) \, dt \geq C \int_{[a,b]^2} (t_1 t_2)^{-d/2} \exp\left(-\frac{\|x\|^2}{C t_1 t_2}\right) dt$$

$$\geq C(b-a)^2 b^{-d} \exp\left(-\frac{M^2}{Ca^2}\right),$$

which shows that Hypothesis H1 is satisfied.

VERIFICATION OF HYPOTHESIS H2.    Fix $x$ and $y$ such that $\|x\| \leq M$ and $\|y\| \leq M$. Let $s, t \in [a, b]^2$ and assume that $s < t$. Clearly,

$$p_{X_t, X_s}(x, y) = p_{X_t - X_s, X_s}(x - y, y)$$

$$= p_{X_t - X_s | X_s = y}(x - y) p_{X_s}(y).$$

Let $z = x - y$ and let $E$ be a Borel subset of $\mathbb{R}^d$. By Proposition 4.13, we have, a.s.,

$$\mathbb{P}\{X_t - X_s \in E | \mathcal{F}_s\} \leq \int_E c \|t - s\|^{-d/2} \exp\left(-\frac{\|z\|^2}{c \|t - s\|}\right) dz$$

for some finite constant $c > 0$. Take the conditional expectation of both sides with respect to $\sigma(X_s)$ to find that

$$p_{X_t - X_s | X_s = y}(z) \leq c \|t - s\|^{-d/2} \exp\left(-\frac{\|z\|^2}{c\|t-s\|}\right).$$

By Lemma 4.12, $p_{X_s}(y)$ is uniformly bounded over $s \in [a,b]^2$ and $x \in \mathbb{R}^d$. Therefore, for all $s, t \in [a,b]^2$ with $s < t$, we have proved that

$$(5.1) \qquad p_{X_t, X_s}(x, y) \leq c' \|t - s\|^{-d/2} \exp\left(-\frac{\|x-y\|^2}{c\|t-s\|}\right)$$

for some finite constant $c' > 0$.

We now assume that $s, t \in [a,b]^2$ with $s \npreceq t$ and $t \npreceq s$. Let $E_1$ and $E_2$ be two Borel subsets of $\mathbb{R}^d$. Using the conditional independence property, we obtain

$$\mathbb{P}\{X_t \in E_1, X_s \in E_2\} = \mathbb{E}[\mathbb{P}\{X_t \in E_1, X_s \in E_2 | \mathcal{F}_{s \wedge t}\}]$$

$$= \mathbb{E}[\mathbb{P}\{X_t \in E_1 | \mathcal{F}_{s \wedge t}\} \mathbb{P}\{X_s \in E_2 | \mathcal{F}_{s \wedge t}\}]$$

$$= \mathbb{E}\left[\int_{E_1} p_{s \wedge t, t}(\cdot, x - X_{s \wedge t}) \, dx \int_{E_2} p_{s \wedge t, s}(\cdot, y - X_{s \wedge t}) \, dy\right].$$

Set $p(c, x) = c^{-d/2} \exp(-\|x\|^2/(2c))$. By Proposition 4.13 and Lemma 4.12, the right-hand side is bounded above by

$$C\mathbb{E}\left[\int_{E_1} dx \, p(c\|t - (s \wedge t)\|, x - X_{s \wedge t}) \int_{E_2} dy \, p(c\|s - (s \wedge t)\|, y - X_{s \wedge t})\right]$$

$$= C \int_{E_1} dx \int_{E_2} dy \int_{\mathbb{R}^d} dz \, p(c\|t - (s \wedge t)\|, x - z) p(c\|s - (s \wedge t)\|, y - z).$$

For $x$ fixed, do the change of variables $u = z - x$ and use the fact that $p(\sigma^2, \cdot)$ is an even function to write the inner integral as a convolution of two Gaussian densities, which is equal to

$$p(c\|t - (s \wedge t)\| + c\|s - (s \wedge t)\|, y - x) \leq Cp(\tilde{c}\|t - s\|, y - x)$$

by the triangle inequality. We conclude that

$$\mathbb{P}\{X_t \in E_1, X_s \in E_2\} \leq C \int_{E_1} dx \int_{E_2} dy \, p(\tilde{c}\|t - s\|, y - x)$$

and, therefore,

$$(5.2) \qquad p_{X_t, X_s}(x, y) \leq C\|t - s\|^{-d/2} \exp\left(-\frac{\|x-y\|^2}{\tilde{c}\|t-s\|}\right).$$

This shows that an estimate like (5.1) holds also when $s \npreceq t$ and $t \npreceq s$.

Using (5.1) and (5.2), we get

$$\int_{[a,b]^2} \int_{[a,b]^2} p_{X_t, X_s}(x, y) \, dt \, ds$$

$$\leq C \int_{[a,b]^2} \int_{[a,b]^2} \|t - s\|^{-d/2} \exp\left(-\frac{\|x - y\|^2}{c\|t - s\|}\right) dt \, ds.$$

Fix $t$ and use the change of variables $u = t - s$ to see that this expression is less than or equal to

$$4C(b - a)^2 \int_{[0, b-a]^2} \|u\|^{-d/2} \exp\left(-\frac{\|x - y\|^2}{c\|u\|}\right) du.$$

Next, use the change of variables $u = \|x - y\|^2 c^{-1} z$, to see that this is less than or equal to

$$c^{d/2 - 2} 4C(b - a)^2 \|x - y\|^{-d+4} \int_{[0, c(b-a)/\|x-y\|^2]^2} \|z\|^{-d/2} \exp(-1/\|z\|) \, dz.$$

When $d \geq 4$, applying Lemma 3.5 with $\beta = d/2 \geq 2$ and $\alpha = 1/2$ in the same way as in the proof of Theorem 3.1 shows that Hypothesis H2 is verified with $k(\cdot) = k_{d-4}(\cdot)$ and $d \geq 4$. When $d \leq 4$, the above expression is bounded and Hypothesis H2 holds with $k(x) \equiv 1$.

VERIFICATION OF HYPOTHESIS H3. Assume $d \geq 4$. Fix $x \in \mathbb{R}^d$ with $\|x\| \leq M$. Use the Gaussian-type lower bound for $p_{s,t}(\omega, x)$, obtained in Theorem 4.16, to see that for all $s \in [a, b]^2$ and for almost all $\omega$,

$$\int_{[b, 2b-a]^2} p_{s,t}(\omega, x) \, dt \geq c \int_{[b, 2b-a]^2} \|t - s\|^{-d/2} \exp\left(-\frac{\|x\|^2}{c\|t - s\|}\right) dt$$

for some finite constant $c > 0$. Now, using the change of variables $z = c(t - s)/\|x\|^2$, we see that this expression is

$$\geq c^{d/2 - 2} \|x\|^{-d+4} \int_{[0, c(b-a)/\|x\|^2]^2} \|z\|^{-d/2} \exp(-1/\|z\|) \, dz.$$

Finally, Lemma 3.5 with $\beta = d/2 \geq 2$ and $\alpha = 1/2$ shows that Hypothesis H3 holds for $d \geq 4$ in the same way as in the proof of Theorem 3.1.

The proof of Theorem 5.1 is complete.  □

5.2. *The case $b \not\equiv 0$.* The aim of this section is to extend Theorem 5.1 to the case $b \not\equiv 0$. Our main tool is Girsanov's theorem for an adapted translation of the Brownian sheet (see [20], Proposition 1.6).

Let $Y = (Y_t, t \in \mathbb{R}_+^2)$ denote the $d$-dimensional adapted continuous process defined on $(\Omega, \mathcal{G}, \mathbb{P})$ satisfying (4.7) with $x_0 = 0$, that is,

$$(5.3) \quad Y_t^i = \sum_{j=1}^d \int_{[0,t]} \sigma_j^i(Y_s)\, dW_s^j + \int_{[0,t]} b^i(Y_s)\, ds, \qquad t \in \mathbb{R}_+^2,\ 1 \le i \le d.$$

We introduce the following condition on the vector $b$:

HYPOTHESIS P3. *For some constant $N$, for all $1 \le i \le d$ and $x \in \mathbb{R}^d$, $|b^i(x)| \le N$.*

Consider the random variable

$$L_t = \exp\left[ -\int_{[0,t]} \sigma^{-1}(Y_s) b(Y_s) \cdot dW_s - \tfrac{1}{2} \int_{[0,t]} \|\sigma^{-1}(Y_s) b(Y_s)\|_{\mathbb{R}^d}^2\, ds \right].$$

We have the following Girsanov theorem.

THEOREM 5.2 ([20], Proposition 1.6). *The random variable $L_t$ is such that $\mathbb{E}[L_t] = 1$. If $\tilde{\mathbb{P}}$ denotes the probability measure on $(\Omega, \mathcal{F})$ defined by*

$$\frac{d\tilde{\mathbb{P}}}{d\mathbb{P}}(\omega) = L_t(\omega),$$

*then $\tilde{W}_t = W_t + \int_{[0,t]} \sigma^{-1}(Y_s) b(Y_s)\, ds$ is a standard Brownian sheet under $\tilde{\mathbb{P}}$.*

Consequently, the law of $Y$ under $\tilde{\mathbb{P}}$ coincides with the law of $X$ under $\mathbb{P}$, where $X = (X_t, t \in \mathbb{R}_+^2)$ is the solution of (5.3) with $b \equiv 0$.

The following result is the extension of Theorem 5.1 to the case $b \not\equiv 0$. It is sufficient to characterize polar sets of $Y$.

COROLLARY 5.3. *Assuming Hypotheses P1–P3, for all $0 < a < b < \infty$ and $\varepsilon, M > 0$, there exists a finite positive constant $K_\varepsilon$ depending on $a$, $b$, $\varepsilon$, $\lambda$, $M$, $N$, $\rho$ and the uniform bounds of the coefficients of $\sigma$ and its derivatives, such that for all compact sets $A \subset \{x \in \mathbb{R}^d : \|x\| < M\}$,*

$$K_\varepsilon^{-1} (\mathrm{Cap}_{d-4}(A))^{1+\varepsilon} \le \mathbb{P}\{\exists\, t \in [a,b]^2 : Y_t \in A\} \le K_\varepsilon (\mathrm{Cap}_{d-4}(A))^{1/1+\varepsilon}.$$

PROOF. Consider the random variable

$$J_t = \exp\left[ -\int_{[0,t]} \sigma^{-1}(X_s) b(X_s) \cdot dW_s + \tfrac{1}{2} \int_{[0,t]} \|\sigma^{-1}(X_s) b(X_s)\|_{\mathbb{R}^d}^2\, ds \right].$$

Fix $0 < a < b < \infty$ and let $G_X = \{\exists\, t \in [a,b]^2 : X_t \in A\}$ and $G_Y = \{\exists\, t \in [a,b]^2 : Y_t \in A\}$. By Theorem 5.2,

$$(5.4) \qquad \mathbb{P}[G_Y] = \mathbb{E}_{\mathbb{P}}[\mathbb{1}_{G_Y}] = \mathbb{E}_{\tilde{\mathbb{P}}}[\mathbb{1}_{G_Y} L_t^{-1}] = \mathbb{E}_{\mathbb{P}}[\mathbb{1}_{G_X} J_t^{-1}].$$

Let $\varepsilon > 0$ and apply Hölder's inequality:

$$\mathbb{P}[G_X] = \mathbb{E}_{\mathbb{P}}[\mathbb{1}_{G_X} J_t^{-1/1+\varepsilon} J_t^{1/1+\varepsilon}] \leq (\mathbb{E}_{\mathbb{P}}[\mathbb{1}_{G_X} J_t^{-1}])^{1/1+\varepsilon} (\mathbb{E}_{\mathbb{P}}[J_t^{1/\varepsilon}])^{\varepsilon/1+\varepsilon}.$$

Rewriting the last inequality we obtain

$$\mathbb{P}[G_Y] \geq (\mathbb{P}[G_X])^{1+\varepsilon} (\mathbb{E}_{\mathbb{P}}[J_t^{1/\varepsilon}])^{-\varepsilon}.$$

Let $r > 0$. By the Cauchy–Schwarz inequality

$$\begin{aligned}
\mathbb{E}_{\mathbb{P}}[J_t^r] \leq \Bigg( &\mathbb{E}_{\mathbb{P}}\Bigg[\exp\Bigg[\int_{[0,t]} (-2)r\sigma^{-1}(X_s)b(X_s) \cdot dW_s \\
&\qquad - \frac{1}{2}\int_{[0,t]} 4r^2 \|\sigma^{-1}(X_s)b(X_s)\|_{\mathbb{R}^d}^2 \, ds\Bigg]\Bigg]\Bigg)^{1/2} \\
&\times \Bigg(\mathbb{E}_{\mathbb{P}}\Bigg[\exp\Bigg[\int_{[0,t]} (2r^2+r)\|\sigma^{-1}(X_s)b(X_s)\|_{\mathbb{R}^d}^2 \, ds\Bigg]\Bigg]\Bigg)^{1/2}.
\end{aligned}$$

The first expectation on the right-hand side equals 1 since it is the expectation of an exponential martingale with bounded quadratic variation (see [9], Chapter 3, Proposition 5.12). By Hypotheses P2 and P3, the second factor is bounded by some positive finite constant. Therefore,

$$\mathbb{E}_{\mathbb{P}}[J_t^{1/\varepsilon}] \leq k_\varepsilon$$

for some constant $k_\varepsilon > 0$. Finally, by Theorem 5.1 there exists a positive finite constant $K$ such that $\mathbb{P}[G_X] \geq K \operatorname{Cap}_{d-4}(A)$, which concludes the proof of the lower bound.

The upper bound is proved along the same lines. Let $\varepsilon > 0$ and apply Hölder's inequality to the right-hand side of (5.4):

$$\mathbb{P}[G_Y] = \mathbb{E}_{\mathbb{P}}[\mathbb{1}_{G_X} J_t^{-1}] \leq (\mathbb{P}[G_X])^{1/1+\varepsilon} (\mathbb{E}_{\mathbb{P}}[J_t^{-(1+\varepsilon/\varepsilon)}])^{\varepsilon/1+\varepsilon}.$$

Let $r > 0$. Again by the Cauchy–Schwarz inequality we find that

$$\begin{aligned}
\mathbb{E}_{\mathbb{P}}[J_t^{-r}] \leq \Bigg( &\mathbb{E}_{\mathbb{P}}\Bigg[\exp\Bigg[\int_{[0,t]} 2r\sigma^{-1}(X_s)b(X_s) \cdot dW_s \\
&\qquad - \frac{1}{2}\int_{[0,t]} 4r^2 \|\sigma^{-1}(X_s)b(X_s)\|_{\mathbb{R}^d}^2 \, ds\Bigg]\Bigg]\Bigg)^{1/2} \\
&\times \Bigg(\mathbb{E}_{\mathbb{P}}\Bigg[\exp\Bigg[\int_{[0,t]} (2r^2-r)\|\sigma^{-1}(X_s)b(X_s)\|_{\mathbb{R}^d}^2 \, ds\Bigg]\Bigg]\Bigg)^{1/2}.
\end{aligned}$$

The first expectation on the right-hand side equals 1 since it is the expectation of an exponential martingale with bounded quadratic variation, as above. By Hypotheses P2 and P3, the second factor is bounded by some positive finite constant. Therefore,

$$\mathbb{E}_{\mathbb{P}}[J_t^{-(1+\varepsilon/\varepsilon)}] \leq k_\varepsilon$$

for some constant $k_\varepsilon > 0$. Finally, by Theorem 5.1, there exists a positive finite positive constant $K$ such that $\mathbb{P}[G_X] \leq K \operatorname{Cap}_{d-4}(A)$, which completes the proof of the Corollary 5.3. $\square$

As a consequence of Corollaries 2.5 and 5.3 we obtain the following analytic criterion for polarity for the process $Y$:

COROLLARY 5.4. *Assume Hypotheses* P1–P3. *Let $E$ be compact subset of $\mathbb{R}^d$. Then $E$ is a polar set for $Y$ if and only if* $\operatorname{Cap}_{d-4}(E) = 0$.

Finally, Theorem 3.3 and Corollary 5.3 give the stochastic codimension and Hausdorff dimension of the range of the process $Y$:

COROLLARY 5.5. *Assuming Hypotheses* P1–P3,

$$\operatorname{codim}\{Y((0, +\infty)^2)\} = (d - 4)^+$$

*and if $d > 4$, then*

$$\dim_{\mathcal{H}}\{Y((0, +\infty)^2)\} = 4 \qquad a.s.$$

5.3. *Critical dimension for hitting points.* For an $N$-parameter $\mathbb{R}^d$-valued Brownian sheet $W = (W_t, t \in \mathbb{R}_+^N)$, Orey and Pruitt [23] showed that for any $x \in \mathbb{R}^d$,

$$\mathbb{P}\{\exists t \in \mathbb{R}_+^N : W_t = x\} = \begin{cases} 1, & \text{if } d < 2N, \\ 0, & \text{if } d \geq 2N. \end{cases}$$

This fact also can be obtained as a consequence of work by Khoshnevisan and Shi [11]. Corollary 5.4 immediately yields an analogous result for the solution $Y$ of (5.3):

COROLLARY 5.6. *Assuming Hypotheses* P1–P3, *for any $x \in \mathbb{R}^d$,*

$$\mathbb{P}\{\exists t \in \mathbb{R}_+^2 : Y_t = x\} > 0 \qquad \text{if and only if } d < 4.$$

**Acknowledgment.** This article is based on E. Nualart's Ph.D. dissertation, written under the supervision of R. C. Dalang.

INSTITUT DE MATHÉMATIQUES
ÉCOLE POLYTECHNIQUE FÉDÉRALE
1015 LAUSANNE
SWITZERLAND
E-MAIL: robert.dalang@epfl.ch
E-MAIL: eulalia.nualart@epfl.ch
URL: http://ima.epfl.ch/prob/membres/dalang.html